# SOME NONASYMPTOTIC RESULTS ON RESAMPLING IN HIGH DIMENSION, I: CONFIDENCE REGIONS[1]

By Sylvain Arlot, Gilles Blanchard[2] and Etienne Roquain

*CNRS ENS, Weierstrass Institut and UPMC University of Paris 6*

We study generalized bootstrap confidence regions for the mean of a random vector whose coordinates have an unknown dependency structure. The random vector is supposed to be either Gaussian or to have a symmetric and bounded distribution. The dimensionality of the vector can possibly be much larger than the number of observations and we focus on a nonasymptotic control of the confidence level, following ideas inspired by recent results in learning theory. We consider two approaches, the first based on a concentration principle (valid for a large class of resampling weights) and the second on a resampled quantile, specifically using Rademacher weights. Several intermediate results established in the approach based on concentration principles are of interest in their own right. We also discuss the question of accuracy when using Monte Carlo approximations of the resampled quantities.

## 1. Introduction.

1.1. *Goals and motivations.* Let $\mathbf{Y} := (\mathbf{Y}^1, \ldots, \mathbf{Y}^n)$ be a sample of $n \geq 2$ i.i.d. observations of an integrable random vector in $\mathbb{R}^K$, with dimensionality $K$ possibly much larger than $n$ and unknown dependency structure of the coordinates. Let $\mu \in \mathbb{R}^K$ denote the common mean of the $\mathbf{Y}^i$. Our goal is to find a nonasymptotic $(1-\alpha)$-confidence region $\mathcal{G}(\mathbf{Y}, 1-\alpha)$ for $\mu$, of the form

(1) $$\mathcal{G}(\mathbf{Y}, 1-\alpha) = \{x \in \mathbb{R}^K | \phi(\overline{\mathbf{Y}} - x) \leq t_\alpha(\mathbf{Y})\},$$

Received November 2007; revised August 2008.
[1]Supported in part by the IST and ICT programs of the European Community, respectively, under the PASCAL (IST-2002-506778) and PASCAL2 (ICT-216886) Networks of Excellence.
[2]Supported in part by the Fraunhofer Institut FIRST, Berlin.
*AMS 2000 subject classifications.* Primary 62G15; secondary 62G09.
*Key words and phrases.* Confidence regions, high-dimensional data, nonasymptotic error control, resampling, cross-validation, concentration inequalities, resampled quantile.

where $\phi : \mathbb{R}^K \to \mathbb{R}$ is a function which is fixed in advance (measuring a kind of distance, e.g., an $\ell_p$-norm for $p \in [1, \infty])$, $\alpha \in (0,1)$, $t_\alpha : (\mathbb{R}^K)^n \to \mathbb{R}$ is a possibly data-dependent threshold and $\overline{\mathbf{Y}} = \frac{1}{n} \sum_{i=1}^n \mathbf{Y}^i \in \mathbb{R}^K$ is the empirical mean of the sample $\mathbf{Y}$.

The point of view developed in the present work focuses on the following goal:

- obtaining *nonasymptotic* results, valid for any fixed $K$ and $n$, with $K$ possibly much larger than the number of observations $n$, while
- avoiding any specific assumption on the dependency structure of the coordinates of $\mathbf{Y}^i$ (although we will consider some general assumptions over the distribution of $\mathbf{Y}$, namely symmetry and boundedness or Gaussianity).

In the Gaussian case, a traditional parametric method based on the direct estimation of the covariance matrix to derive a confidence region would not be appropriate in the situation where $K \gg n$, unless the covariance matrix is assumed to belong to some parametric model of lower dimension, which we explicitly do not want to posit here. In this sense, the approach followed here is closer in spirit to nonparametric or semiparametric statistics.

This point of view is motivated by some practical applications, especially neuroimaging [8, 18, 26]. In a magnetoencephalography (MEG) experiment, each observation $\mathbf{Y}^i$ is a two- or three-dimensional brain activity map, obtained as a difference between brain activities in the presence or absence of some stimulation. The activity map is typically composed of about 15,000 points; the data can also be a time series of length between 50 and 1000 such maps. The dimensionality $K$ can thus range from $10^4$ to $10^7$. Such observations are repeated from $n = 15$ up to 4000 times, but this upper bound is seldom attained [32]; in typical cases, one has $n \leq 100 \ll K$. In such data, there are strong dependencies between locations (the 15,000 points are obtained by pre-processing data from 150 sensors) and these dependencies are spatially highly nonhomogeneous, as noted by [26]. Moreover, there may be long-distance correlations, for example, depending on neural connections inside the brain, so that a simple parametric model of the dependency structure is generally not adequate. Another motivating example is given by microarray data [14], where it is common to observe samples of limited size (say, less than 100) of a vector in high dimension (say, more than 20,000, each dimension corresponding to a specific gene) and where the dependency structure may be quite arbitrary.

1.2. *Two approaches to our goal.* The ideal threshold $t_\alpha$ in (1) is obviously the $(1-\alpha)$ quantile of the distribution of $\phi(\overline{\mathbf{Y}} - \mu)$. However, this quantity depends on the unknown dependency structure of the coordinates of $\mathbf{Y}^i$ and is therefore itself unknown.

The approach studied in this work is to use a (generalized) resampling scheme in order to estimate $t_\alpha$. The heuristics of the resampling method (introduced in [11], generalized to exchangeable weighted bootstrap by [23, 28]) is that the distribution of the unobservable variable $\overline{\mathbf{Y}} - \mu$ is "mimicked" by the distribution, conditionally on $\mathbf{Y}$, of the resampled empirical mean of the centered data. This last quantity is an observable variable and we denote it as follows:

$$(2) \quad \overline{\mathbf{Y}}^{\langle W - \overline{W} \rangle} := \frac{1}{n} \sum_{i=1}^{n} (W_i - \overline{W}) \mathbf{Y}^i = \frac{1}{n} \sum_{i=1}^{n} W_i (\mathbf{Y}^i - \overline{\mathbf{Y}}) = (\overline{\mathbf{Y} - \overline{\mathbf{Y}}})^{\langle W \rangle},$$

where $(W_i)_{1 \leq i \leq n}$ are real random variables independent of $\mathbf{Y}$, called the *resampling weights*, and $\overline{W} = n^{-1} \sum_{i=1}^{n} W_i$. We emphasize that the weight family $(W_i)_{1 \leq i \leq n}$ itself *need not be independent*.

In Section 2.4, we define in more detail several specific resampling weights inspired both by traditional resampling methods [23, 28] and by recent statistical learning theory. Let us give two typical examples reflecting these two sources:

- *Efron's bootstrap weights*. $W$ is a multinomial random vector with parameters $(n; n^{-1}, \ldots, n^{-1})$. This is the standard bootstrap.
- *Rademacher weights*. $W_i$ are i.i.d. Rademacher variables, that is, $W_i \in \{-1, 1\}$ with equal probabilities. They are closely related to symmetrization techniques in learning theory.

It is useful to observe at this point that, to the extent that we only consider resampled data after empirical centering, shifting all weights by the same (but possibly random) offset $C > 0$ does not change the resampled quantity introduced in (2). Hence, to reconcile the intuition of traditional resampling with what could possibly appear as unfamiliar weights, one could always assume that the weights are translated to enforce (for example) weight positivity or the condition $n^{-1} \sum_{i=1}^{n} W_i = 1$ (although, of course, in general, both conditions cannot be ensured at the same time simply by translation). For example, Rademacher weights can be interpreted as a resampling scheme where each $\mathbf{Y}^i$ is independently discarded or "doubled" with equal probability.

Following the general resampling idea, we investigate two distinct approaches in order to obtain nonasymptotic confidence regions:

- Approach 1 ("concentration approach," developed in Section 2).

    The expectations of $\phi(\overline{\mathbf{Y}} - \mu)$ and $\phi(\overline{\mathbf{Y}}^{\langle W - \overline{W} \rangle})$ can be precisely compared and the processes $\phi(\overline{\mathbf{Y}} - \mu)$ and $\mathbb{E}_W[\phi(\overline{\mathbf{Y}}^{\langle W - \overline{W} \rangle})]$ concentrate well around their respective expectations, where $\mathbb{E}_W$ denotes the expectation operator with respect to the distribution of $W$ (i.e., conditionally on $\mathbf{Y}$).

- Approach 2 ("direct quantile approach," developed in Section 3).

  The $1 - \alpha$ quantile of the distribution of $\phi(\overline{\mathbf{Y}}^{\langle W - \overline{W} \rangle})$ conditionally on $\mathbf{Y}$ is close to the $1 - \alpha$ quantile of $\phi(\overline{\mathbf{Y}} - \mu)$.

Regarding the second approach, we will restrict ourselves specifically to Rademacher weights in our analysis and rely heavily on a symmetrization principle.

1.3. *Relation to previous work.* Using resampling to construct confidence regions is a vast field of study in statistics (see, e.g., [4, 9, 11, 15, 16, 27]). Available results are, however, mostly asymptotic, based on the celebrated fact that the bootstrap process is asymptotically close to the original empirical process [31]. Because we focus on a nonasymptotic viewpoint, this asymptotic approach is not adapted to the goals we have fixed. Note, also, that the nonasymptotic viewpoint can be used as a basis for an asymptotic analysis in the situation where the dimension $K$ grows with $n$, a setting which is typically not covered by standard asymptotics.

The "concentration approach" mentioned in the previous subsection is inspired by recent results coming from learning theory and relates, in particular, to the notion of Rademacher complexity [20]. This notion has been extended in the recent work of Fromont [13] to more general resampling schemes and this latter work has had a strong influence on the present paper.

On the other hand, what we called the "quantile approach" in the previous subsection is strongly related to exact randomization tests (which are based on an invariance of the null distribution under a given transformation; the underlying idea can be traced back to Fisher's permutation test [12]). Namely, we will only consider symmetric distributions; this is a specific instance of an invariance with respect to a transformation and will allow us to make use of distribution-preserving randomization via sign reversal. Here, the main difference with traditional exact randomization tests is that, since our goal is to derive a confidence region, the vector of the means is unknown and, therefore, so is the exact invariant transformation. Our contribution to this point is essentially to show that the true vector of the means can be replaced by the empirical one in the randomization, at the cost of additional terms of smaller order in the threshold thus obtained. To our knowledge, this gives the first nonasymptotic approximation result on resampled quantiles with an unknown distribution mean.

Finally, we contrast the setting studied here with a strand of research studying adaptive confidence regions (in a majority of cases, $\ell_2$-balls) in nonparametric Gaussian regression. A seminal paper on this topic is [22] and recent work includes [17, 21, 29] (from an asymptotic point of view) and [3, 5, 6, 19] (which present nonasymptotic results). Related to this setting and ours is [10], where adaptive tests for zero mean are developed for symmetric distributions, using randomization by sign reversal. The setting considered in these papers is that of regression on a fixed design in high dimension (or in the Gaussian sequence model) with one observation per point and i.i.d. noise. This corresponds (in our notation) to $n=1$, while the $K$ coordinates are assumed independent. Despite some similarities, the problem considered here is of a different nature: in the aforementioned works, the focus is on adaptivity with respect to some properties of the true mean vector, concretized by a family of models (e.g., linear subspaces or Besov balls in the Gaussian sequence setting); usually, an adaptive estimator performing implicit or explicit model selection relative to this collection is studied and a crucial question for obtaining confidence regions is that of empirically estimating the bias of this estimator when the noise dependence structure is known. In the present paper, we do not consider the problem of model selection, but the focus is on evaluating the estimation error under an unknown noise dependence structure (for the "naive" unbiased estimator given by the empirical mean).

1.4. *Notation.* We first introduce some notation that will be useful throughout the paper.

- A boldface letter indicates a matrix. This will almost exclusively concern the $K \times n$ data matrix $\mathbf{Y}$. A superscript index such as $\mathbf{Y}^i$ indicates the $i$th column of a matrix.
- If $\mu \in \mathbb{R}^K$, $\mathbf{Y} - \mu$ is the matrix obtained by subtracting $\mu$ from each (column) vector of $\mathbf{Y}$. Similarly, for a vector $W \in \mathbb{R}^n$ and $c \in \mathbb{R}$, we denote $W - c := (W_i - c)_{1 \leq i \leq n} \in \mathbb{R}^n$.
- If $X$ is a random variable, then $\mathcal{D}(X)$ is its distribution and $\text{Var}(X)$ is its variance. We use the notation $X \sim Y$ to indicate that $X$ and $Y$ have the same distribution. Moreover, the support of $\mathcal{D}(X)$ is denoted by $\text{supp}\,\mathcal{D}(X)$.
- We denote by $\mathbb{E}_W[\cdot]$ the expectation operator over the distribution of the weight vector $W$ only, that is, conditional on $\mathbf{Y}$. We use a similar notation, $\mathbb{P}_W$, for the corresponding probability operator and $\mathbb{E}_\mathbf{Y}, \mathbb{P}_\mathbf{Y}$ for the same operations conditional on $W$. Since $\mathbf{Y}$ and $W$ are always assumed to be independent, the operators $\mathbb{E}_W$ and $\mathbb{E}_\mathbf{Y}$ commute by Fubini's theorem.
- The vector $\sigma = (\sigma_k)_{1 \leq k \leq K}$ is the vector of the standard deviations of the data: $\forall k, 1 \leq k \leq K$, $\sigma_k := \text{Var}^{1/2}(\mathbf{Y}_k^1)$.
- $\overline{\Phi}$ is the standard Gaussian upper tail function: if $X \sim \mathcal{N}(0,1)$, $\forall x \in \mathbb{R}$, $\overline{\Phi}(x) := \mathbb{P}(X \geq x)$.
- We define the mean of the weight vector $\overline{W} := \frac{1}{n} \sum_{i=1}^n W_i$, the empirical mean vector $\overline{\mathbf{Y}} := \frac{1}{n} \sum_{i=1}^n \mathbf{Y}^i$ and the resampled empirical mean vector $\overline{\mathbf{Y}}^{\langle W \rangle} := \frac{1}{n} \sum_{i=1}^n W_i \mathbf{Y}^i$.

- We use the operator $|\cdot|$ to denote the cardinality of a set.
- For two positive sequences $(u_n)_n$ and $(v_n)_n$, we write $u_n = \Theta(v_n)$ when $(u_n v_n^{-1})_n$ stays bounded away from zero and $+\infty$.

Several properties may be assumed for the function $\phi : \mathbb{R}^K \to \mathbb{R}$ used to define confidence regions of the form (1):

- Subadditivity: $\forall x, x' \in \mathbb{R}^K, \phi(x + x') \leq \phi(x) + \phi(x')$.
- Positive homogeneity: $\forall x \in \mathbb{R}^K, \forall \lambda \in \mathbb{R}^+, \phi(\lambda x) = \lambda \phi(x)$.
- Boundedness by the $\ell_p$-norm, $p \in [1, \infty]$: $\forall x \in \mathbb{R}^K$, $|\phi(x)| \leq \|x\|_p$, where $\|x\|_p$ is equal to $(\sum_{k=1}^{K} |x_k|^p)^{1/p}$ if $p < \infty$ and $\max_k\{|x_k|\}$ for $p = +\infty$. Note, also, that all of the results in this paper are still valid with any normalization of the $\ell_p$-norm [in particular, it can be taken to be equal to $(K^{-1} \sum_{k=1}^{K} |x_k|^p)^{1/p}$ so that the $\ell_p$-norm of a vector with equal coordinates does not depend on the dimensionality $K$].

Finally, we introduce the following possible assumptions on the generating distribution of $\mathbf{Y}$:

(GA) The Gaussian assumption: the $\mathbf{Y}^i$ are Gaussian vectors.
(SA) The symmetric assumption: the $\mathbf{Y}^i$ are symmetric with respect to $\mu$, that is, $(\mathbf{Y}^i - \mu) \sim (\mu - \mathbf{Y}^i)$.
(BA) $(p, M)$ the boundedness assumption: $\|\mathbf{Y}^i - \mu\|_p \leq M$ a.s.

In this paper, we primarily focus on the Gaussian framework (GA), where the corresponding results will be more accurate. In the sequel, when considering (GA) and the assumption that $\phi$ is bounded by the $\ell_p$-norm for some $p \geq 1$, we will additionally always assume that we know some upper bound on the $\ell_p$-norm of $\sigma$. The question of finding an upper bound for $\|\sigma\|_p$ based on the data is discussed in Section 4.1.

**2. Confidence region using concentration.**

2.1. *Main result.* We consider here a general *resampling weight vector* $W$, that is, an $\mathbb{R}^n$-valued random vector $W = (W_i)_{1 \leq i \leq n}$ independent of $\mathbf{Y}$ and satisfying the following properties: for all $i \in \{1, \ldots, n\}$, $\mathbb{E}[W_i^2] < \infty$ and $n^{-1} \sum_{i=1}^{n} \mathbb{E}|W_i - \overline{W}| > 0$.

In this section, we will mainly consider an *exchangeable resampling weight vector*, that is, a resampling weight vector $W$ such that $(W_i)_{1 \leq i \leq n}$ has an exchangeable distribution (in other words, it is invariant under any permutation of the indices). Several examples of exchangeable resampling weight vectors are given in Section 2.4, where we also address the question of how to choose between different possible distributions of $W$. An extension of our results to nonexchangeable weight vectors is proposed in Section 2.5.1.

TABLE 1
*Resampling constants for some classical resampling weight vectors*

| | |
|---|---|
| Efron | $2(1-\frac{1}{n})^n = A_W \leq B_W \leq \sqrt{\frac{n-1}{n}}$, $C_W = 1$ |
| Efron, $n \to +\infty$ | $\frac{2}{e} = A_W \leq B_W \leq 1 = C_W$ |
| Rademacher | $1 - \frac{1}{n} = A_W \leq B_W \leq \sqrt{1-\frac{1}{n}}$, $C_W = 1 \leq D_W \leq 1 + \frac{1}{\sqrt{n}}$ |
| Rademacher, $n \to +\infty$ | $A_W = B_W = C_W = D_W = 1$ |
| rho($q$) | $A_W = 2(1 - \frac{q}{n})$, $B_W = \sqrt{\frac{n}{q} - 1}$ |
| | $C_W = \sqrt{\frac{n}{n-1}}\sqrt{\frac{n}{q} - 1}$, $D_W = \frac{n}{2q} + |1 - \frac{n}{2q}|$ |
| rho($n/2$) | $A_W = B_W = D_W = 1$, $C_W = \sqrt{\frac{n}{n-1}}$ |
| Leave-one-out | $\frac{2}{n} = A_W \leq B_W = \frac{1}{\sqrt{n-1}}$, $C_W = \frac{\sqrt{n}}{n-1}$, $D_W = 1$ |
| Regular $V$-fcv | $A_W = \frac{2}{V} \leq B_W = \frac{1}{\sqrt{V-1}}$, $C_W = \sqrt{n}(V-1)^{-1}$, $D_W = 1$ |

Four constants that depend only on the distribution of $W$ appear in the results below (the fourth one is defined only for a particular class of weights). They are defined as follows and computed for classical resamplings in Table 1:

$$A_W := \mathbb{E}|W_1 - \overline{W}|; \tag{3}$$

$$B_W := \mathbb{E}\left[\left(\frac{1}{n}\sum_{i=1}^n (W_i - \overline{W})^2\right)^{1/2}\right]; \tag{4}$$

$$C_W := \left(\frac{n}{n-1}\mathbb{E}[(W_1 - \overline{W})^2]\right)^{1/2}; \tag{5}$$

$$D_W := a + \mathbb{E}|\overline{W} - x_0| \tag{6}$$

if, for all $i, |W_i - x_0| = a$ a.s. (with $a > 0, x_0 \in \mathbb{R}$).

Note that these quantities are positive for an exchangeable resampling weight vector $W$ and satisfy

$$0 < A_W \leq B_W \leq C_W \sqrt{1 - 1/n}.$$

Moreover, if the weights are i.i.d., we have $C_W = \mathrm{Var}(W_1)^{1/2}$. We can now state the main result of this section.

THEOREM 2.1. *Fix $\alpha \in (0,1)$ and $p \in [1, \infty]$. Let $\phi: \mathbb{R}^K \to \mathbb{R}$ be any function which is subadditive, positive homogeneous and bounded by the $\ell_p$-norm, and let $W$ be an exchangeable resampling weight vector.*

1. *If* $\mathbf{Y}$ *satisfies* (GA), *then*

$$
(7) \qquad \phi(\overline{\mathbf{Y}} - \mu) < \frac{\mathbb{E}_W[\phi(\overline{\mathbf{Y}}^{\langle W - \overline{W} \rangle})]}{B_W} + \|\sigma\|_p \overline{\Phi}^{-1}(\alpha/2) \left[ \frac{C_W}{nB_W} + \frac{1}{\sqrt{n}} \right]
$$

holds with probability at least $1 - \alpha$. The same bound holds for the lower deviations, that is, with inequality (7) reversed and the additive term replaced by its opposite.

2. *If* $\mathbf{Y}$ *satisfies* (SA) *and* (BA) $(p, M)$ *for some* $M > 0$, *then*

$$
(8) \qquad \phi(\overline{\mathbf{Y}} - \mu) < \frac{\mathbb{E}_W[\phi(\overline{\mathbf{Y}}^{\langle W - \overline{W} \rangle})]}{A_W} + \frac{2M}{\sqrt{n}} \sqrt{\log(1/\alpha)}
$$

holds with probability at least $1 - \alpha$. If, moreover, the weight vector satisfies the assumption of (6), then

$$
(9) \qquad \phi(\overline{\mathbf{Y}} - \mu) > \frac{\mathbb{E}_W[\phi(\overline{\mathbf{Y}}^{\langle W - \overline{W} \rangle})]}{D_W} - \frac{M}{\sqrt{n}} \sqrt{1 + \frac{A_W^2}{D_W^2}} \sqrt{2 \log(1/\alpha)}
$$

holds with probability at least $1 - \alpha$.

Inequalities (7) and (8) give regions of the form (1) that are confidence regions of level at least $1 - \alpha$. They require knowledge of some upper bound on $\|\sigma\|_p$ (resp., $M$) or a good estimate of it. We address this question in Section 4.1.

In order to obtain some insight into these bounds, it is useful to compare them with an elementary inequality. In the Gaussian case, it is true for each coordinate $k, 1 \leq k \leq K$, that the following inequality holds with probability $1 - \alpha$: $|\overline{\mathbf{Y}}_k - \mu_k| < \frac{\sigma_k}{\sqrt{n}} \overline{\Phi}^{-1}(\alpha/2)$. By applying a simple union bound over the coordinates and using the fact that $\phi$ is positive homogeneous and bounded by the $\ell_p$-norm, we conclude that the following inequality holds with probability at least $1 - \alpha$:

$$
(10) \qquad \phi(\overline{\mathbf{Y}} - \mu) < \frac{\|\sigma\|_p}{\sqrt{n}} \overline{\Phi}^{-1}\left(\frac{\alpha}{2K}\right) =: t_{\text{Bonf}}(\alpha),
$$

which is a minor variation on the well-known Bonferroni bound. By comparison, the main term in the remainder part of (7) takes a similar form, but with $K$ replaced by 1: the remainder term is *dimension-independent*. Naturally, the "dimension complexity" has not disappeared, but will be taken into account in the main resampled term instead. When $K$ is large, the bound (7) can improve on the Bonferroni threshold if there are strong dependencies between the coordinates, resulting in a significantly smaller resampling term.

By way of illustration, consider an extreme example where all pairwise coordinate correlations are exactly 1, that is, the random vector $\mathbf{Y}$ is made of $K$ copies of the same random variable so that there is, in fact, no dimension complexity. Take $\phi(X) = \sup_i X_i$ (corresponding to a uniform one-sided confidence bound for the mean components). Then the resampled quantity in (7) is equal to zero and the obtained bound is close to optimal (up to the two following points: the level is divided by a factor of 2 and there is an additional term of order $\frac{1}{n}$). By comparison, the Bonferroni bound divides the level by a factor of $K$, resulting in a significantly worse threshold. In passing, note that this example illustrates that the order $n^{-1/2}$ of the remainder term cannot be improved.

If we now interpret the bound (7) from an asymptotic point of view [with $K(n)$ depending on $n$ and $\|\sigma\|_p = \Theta(1)$], then the rate of convergence to zero cannot be faster than $n^{-1/2}$ (which corresponds to the standard parametric rate when $K$ is fixed), but it can be potentially slower, for example, if $K$ increases exponentially with $n$. In the latter case, the rate of convergence of the Bonferroni threshold is always strictly slower than $n^{-1/2}$. In general, as far as the order in $n$ is concerned, the resampled threshold converges at least as fast as Bonferroni's, but whether it is strictly faster depends once again on the coordinate dependency structure.

However, if the coordinates are only "weakly dependent," then the threshold (7) can be more conservative than Bonferroni's by a multiplicative factor, while the Bonferroni threshold can sometimes be essentially optimal (e.g., with $\phi = \|\cdot\|_\infty$, all of the coordinates independent and with small $\alpha$). This motivates the next result, where we assume, more generally, that an alternate analysis of the problem can lead to deriving a *deterministic* threshold $t_\alpha$ such that $\mathbb{P}(\phi(\overline{\mathbf{Y}} - \mu) > t_\alpha) \leq \alpha$. In this case, we would ideally like to take the "best of two approaches" and consider the minimum of $t_\alpha$ and the resampling-based thresholds considered above. In the Gaussian case, the following proposition establishes that we can combine the concentration threshold corresponding to (7) with $t_\alpha$ to obtain a threshold that is very close to the minimum of the two.

PROPOSITION 2.2. *Fix $\alpha, \delta \in (0,1)$, $p \in [1, \infty]$ and take $\phi$ and $W$ as in Theorem 2.1. Suppose that $\mathbf{Y}$ satisfies* (GA) *and that $t_{\alpha(1-\delta)}$ is a real number such that $\mathbb{P}(\phi(\overline{\mathbf{Y}} - \mu) > t_{\alpha(1-\delta)}) \leq \alpha(1-\delta)$. Then, with probability at least $1 - \alpha$, $\phi(\overline{\mathbf{Y}} - \mu)$ is less than or equal to the minimum of $t_{\alpha(1-\delta)}$ and*

$$(11) \quad \frac{\mathbb{E}_W[\phi(\overline{\mathbf{Y}}^{\langle W - \overline{W}\rangle})]}{B_W} + \frac{\|\sigma\|_p}{\sqrt{n}}\overline{\Phi}^{-1}\left(\frac{\alpha(1-\delta)}{2}\right) + \frac{\|\sigma\|_p C_W}{nB_W}\overline{\Phi}^{-1}\left(\frac{\alpha\delta}{2}\right).$$

The important point to note in Proposition 2.2 is that, since the last term of (11) becomes negligible with respect to the rest when $n$ grows large, we can choose $\delta$ to be quite small [typically $\delta = \Theta(1/n)$] and obtain a threshold very close to the minimum of $t_\alpha$ and the threshold corresponding to (7). Therefore, this result is more subtle than just considering the minimum of two thresholds each taken at level $1 - \frac{\alpha}{2}$, as would be obtained by a direct union bound.

The proof of Theorem 2.1 involves results which are of interest in their own right: the comparison between the expectations of the two processes $\mathbb{E}_W[\phi(\overline{\mathbf{Y}}^{\langle W - \overline{W} \rangle})]$ and $\phi(\overline{\mathbf{Y}} - \mu)$ and the concentration of these processes around their means. These two issues are, respectively, examined in the two next Sections 2.2 and 2.3. In Section 2.4, we provide some elements for an appropriate choice of resampling weight vectors among several classical examples. The final Section 2.5 tackles the practical issue of computation time.

2.2. *Comparison in expectation.* In this section, we compare $\mathbb{E}[\phi(\overline{\mathbf{Y}}^{\langle W - \overline{W} \rangle})]$ and $\mathbb{E}[\phi(\overline{\mathbf{Y}} - \mu)]$. We note that these expectations exist in the Gaussian (GA) and the bounded (BA) cases, provided that $\phi$ is measurable and bounded by an $\ell_p$-norm. Otherwise (in particular, in Propositions 2.3 and 2.4), we assume that these expectations exist.

In the Gaussian case, these quantities are equal up to a factor that depends only on the distribution of $W$.

PROPOSITION 2.3. *Let $\mathbf{Y}$ be a sample satisfying* (GA) *and let $W$ be a resampling weight vector. Then, for any measurable positive homogeneous function $\phi : \mathbb{R}^K \to \mathbb{R}$, we have the following equality:*

$$(12) \qquad B_W \mathbb{E}[\phi(\overline{\mathbf{Y}} - \mu)] = \mathbb{E}[\phi(\overline{\mathbf{Y}}^{\langle W - \overline{W} \rangle})].$$

*If the weights are such that $\sum_{i=1}^n (W_i - \overline{W})^2 = n$, then the above equality holds for any function $\phi$ (and $B_W = 1$).*

For some classical weights, we give bounds or exact expressions for $B_W$ in Table 1. In general, we can compute the value of $B_W$ by simulation. Note that in a non-Gaussian framework, the constant $B_W$ is still of interest, in an asymptotic sense: Theorem 3.6.13 in [31] uses the limit of $B_W$ when $n$ goes to infinity as a normalizing constant.

When the sample is only assumed to have a symmetric distribution, we obtain the following inequalities.

PROPOSITION 2.4. *Let $\mathbf{Y}$ be a sample satisfying* (SA), *$W$ an exchangeable resampling weight vector and $\phi : \mathbb{R}^K \to \mathbb{R}$ any subadditive, positive homogeneous function.*

(i) *We have the following general lower bound:*

$$(13) \qquad A_W \mathbb{E}[\phi(\overline{\mathbf{Y}} - \mu)] \leq \mathbb{E}[\phi(\overline{\mathbf{Y}}^{\langle W - \overline{W} \rangle})].$$

(ii) *If the weight vector satisfies the assumption of (6), then we have the following upper bound:*

$$D_W \mathbb{E}[\phi(\overline{\mathbf{Y}} - \mu)] \geq \mathbb{E}[\phi(\overline{\mathbf{Y}}^{\langle W - \overline{W} \rangle})]. \tag{14}$$

The bounds (13) and (14) are tight (i.e., $A_W/D_W \to 1$ as $n \to \infty$) for some classical weights; see Table 1. When $\mathbf{Y}$ is not assumed to have a symmetric distribution and $\overline{W} = 1$ a.s., Proposition 2 of [13] shows that (13) holds with $A_W$ replaced by $\mathbb{E}(W_1 - \overline{W})_+$. Therefore, assumption (SA) allows us to get a tighter result (e.g., twice as sharp with Efron or Rademacher weights). It can be shown (see [1], Chapter 9) that this factor of 2 is unavoidable in general for a fixed $n$ when (SA) is not satisfied, although it is unnecessary when $n$ goes to infinity. We conjecture that an inequality close to (13) holds under an assumption less restrictive than (SA) (e.g., concerning an appropriate measure of skewness of the distribution of $\mathbf{Y}^1$).

2.3. *Concentration around the expectation.* In this section, we present concentration results for the two processes $\phi(\overline{\mathbf{Y}} - \mu)$ and $\mathbb{E}_W[\phi(\overline{\mathbf{Y}}^{\langle W - \overline{W} \rangle})]$.

PROPOSITION 2.5. *Let $p \in [1, \infty]$, $\mathbf{Y}$ be a sample satisfying* (GA) *and $\phi : \mathbb{R}^K \to \mathbb{R}$ be any subadditive function, bounded by the $\ell_p$-norm.*

(i) *For all $\alpha \in (0, 1)$, with probability at least $1 - \alpha$, we have*

$$\phi(\overline{\mathbf{Y}} - \mu) < \mathbb{E}[\phi(\overline{\mathbf{Y}} - \mu)] + \frac{\|\sigma\|_p \overline{\Phi}^{-1}(\alpha/2)}{\sqrt{n}} \tag{15}$$

*and the same bound holds for the corresponding lower deviations.*

(ii) *Let $W$ be an exchangeable resampling weight vector. Then, for all $\alpha \in (0, 1)$, with probability at least $1 - \alpha$, we have*

$$\mathbb{E}_W[\phi(\overline{\mathbf{Y}}^{\langle W - \overline{W} \rangle})] < \mathbb{E}[\phi(\overline{\mathbf{Y}}^{\langle W - \overline{W} \rangle})] + \frac{\|\sigma\|_p C_W \overline{\Phi}^{-1}(\alpha/2)}{n} \tag{16}$$

*and the same bound holds for the corresponding lower deviations.*

The bound (15) with a remainder in $n^{-1/2}$ is classical; this order in $n$ cannot be improved, as seen, for example, by taking $K = 1$ and $\phi$ to be the identity function. The bound (16) is more interesting because it illustrates one of the key properties of resampling, the "stabilization effect": the resampled expectation concentrates much faster to its expectation than the original quantity. This effect is known and has been studied asymptotically (in fixed dimension) using Edgeworth expansions (see [15]); here, we demonstrate its validity nonasymptotically in a specific case (see also Section 4.2 below for additional discussion).

In the bounded case, the next proposition is a minor variation of a result by Fromont. It is a consequence of McDiarmid's inequality [25]; we refer the reader to [13] (Proposition 1) for a proof.

PROPOSITION 2.6. *Let $p \in [1, \infty]$, $M > 0$, $\mathbf{Y}$ be a sample satisfying* (BA) *$(p, M)$ and $\phi : \mathbb{R}^K \to \mathbb{R}$ be any subadditive function bounded by the $\ell_p$-norm.*

(i) *For all $\alpha \in (0, 1)$, with probability at least $1 - \alpha$, we have*

$$\phi(\overline{\mathbf{Y}} - \mu) < \mathbb{E}[\phi(\overline{\mathbf{Y}} - \mu)] + \frac{M}{\sqrt{n}} \sqrt{\log(1/\alpha)} \tag{17}$$

*and the same bound holds for the corresponding lower deviations.*

(ii) *Let $W$ be an exchangeable resampling weight vector. Then, for all $\alpha \in (0, 1)$, with probability at least $1 - \alpha$, we have*

$$\mathbb{E}_W[\phi(\overline{\mathbf{Y}}^{\langle W - \overline{W}\rangle})] < \mathbb{E}[\phi(\overline{\mathbf{Y}}^{\langle W - \overline{W}\rangle})] + \frac{A_W M}{\sqrt{n}} \sqrt{\log(1/\alpha)} \tag{18}$$

*and the same bound holds for the corresponding lower deviations.*

2.4. *Resampling weight vectors.* In this section, we consider the question of choosing an appropriate exchangeable resampling weight vector $W$ when using Theorem 2.1 or Corollary 2.2. We define the following resampling weight vectors:

1. *Rademacher.* $W_i$ i.i.d. Rademacher variables, that is, $W_i \in \{-1, 1\}$ with equal probabilities.
2. *Efron* (*Efron's bootstrap weights*). $W$ has a multinomial distribution with parameters $(n; n^{-1}, \ldots, n^{-1})$.
3. *Random hold-out(q)* [rho(q) for short], $q \in \{1, \ldots, n\}$. $W_i = \frac{n}{q} \mathbb{1}_{i \in I}$, where $I$ is uniformly distributed on subsets of $\{1, \ldots, n\}$ of cardinality $q$. These weights may also be called *cross-validation weights* or *leave-$(n-q)$-out weights*. A classical choice is $q = n/2$ (assuming $n$ is even). When $q = n - 1$, these weights are called *leave-one-out* weights. Note that this resampling scheme is a particular case of subsampling.

As noted in the Introduction, the first example is common in learning theory, while the second is classical in the framework of the resampling literature [23, 28]. Random hold-out weights have the particular quality of being related to both: they are nonnegative, satisfy $\sum_i W_i = n$ a.s. and originate with a data-splitting idea (choosing $I$ amounts to choose a subsample) upon which the cross-validation idea has been built. This analogy motivates the "$V$-fold cross-validation weights" (defined in Section 2.5), in order to reduce the computational complexity of the procedures proposed here.

For these classical weights, exact or approximate values for the quantities $A_W$, $B_W$, $C_W$ and $D_W$ [defined by (3) to (6)] can easily be derived (see

Table 1). Proofs are given in Section 5.3, where several other weights are considered. Now, to use Theorem 2.1 or Corollary 2.2, we have to choose a particular resampling weight vector. In the Gaussian case, we propose the following accuracy and complexity criteria:

- First, relation (7) suggests that the quantity $C_W B_W^{-1}$ can be proposed as an *accuracy* index for $W$. Namely, this index enters directly into the deviation term of the upper bound (while we know from Proposition 2.3 that the expectation term is exact) so that the smaller this index is, the sharper the bound.
- Second, an upper bound on the computational burden of exactly computing the resampling quantity is given by the cardinality of the support of $\mathcal{D}(W)$, thus providing a *complexity* index.

These two criteria are estimated in Table 2 for classical weights. For any exchangeable weight vector $W$, we have $C_W B_W^{-1} \geq [n/(n-1)]^{1/2}$ and the cardinality of the support of $\mathcal{D}(W)$ is larger than $n$. Therefore, the *leave-one-out weights* satisfy the best accuracy/complexity trade-off among exchangeable weights.

2.5. *Practical computation of the thresholds.* In practice, the exact computation of the resampling quantity $\mathbb{E}_W[\phi(\overline{\mathbf{Y}}^{\langle W-\overline{W}\rangle})]$ can still be too complex for the weights defined above. In this section, we consider two possible ways to address this issue. First, it is possible to use nonexchangeable weights with a lower complexity index and for which the exact computation is tractable. Alternatively, we propose to use a Monte Carlo approximation, as is often done in practice to compute resampled quantities. In both cases, the thresholds have to be made slightly larger in order to keep a rigourous nonasymptotic control on the level. This is detailed in the two paragraphs below.

TABLE 2
*Choice of the resampling weight vectors: accuracy/complexity trade-off*

| **Resampling** | $C_W B_W^{-1}$ **(accuracy)** | $\|\mathrm{supp}\,\mathcal{D}(W)\|$ **(complexity)** |
|---|---|---|
| Efron | $\leq \frac{1}{2}(1-\frac{1}{n})^{-n} \xrightarrow[n\to\infty]{} \frac{e}{2}$ | $\binom{2n-1}{n-1} = \Theta(n^{-1/2} 4^n)$ |
| Rademacher | $\leq n/(n-1) \xrightarrow[n\to\infty]{} 1$ | $2^n$ |
| rho($n/2$) | $= \sqrt{\frac{n}{n-1}} \xrightarrow[n\to\infty]{} 1$ | $\binom{n}{n/2} = \Theta(n^{-1/2} 2^n)$ |
| Leave-one-out | $= \sqrt{\frac{n}{n-1}} \xrightarrow[n\to\infty]{} 1$ | $n$ |
| Regular $V$-fcv | $= \sqrt{\frac{n}{V-1}}$ | $V$ |

2.5.1. *V-fold cross-validation weights.* In order to reduce the computation complexity, we can use "piecewise exchangeable" weights: consider a regular partition $(B_j)_{1 \leq j \leq V}$ of $\{1, \ldots, n\}$ (where $V \in \{2, \ldots, n\}$ and $V$ divides $n$) and define the weights $W_i = \frac{V}{V-1} \mathbb{1}_{i \notin B_J}$ with $J$ uniformly distributed on $\{1, \ldots, V\}$. These weights are called the (*regular*) *V-fold cross-validation* weights (*V*-fcv for short).

By applying our previous results to the process $(\widetilde{\mathbf{Y}}^j)_{1 \leq j \leq V}$, where $\widetilde{\mathbf{Y}}^j := \frac{V}{n} \sum_{i \in B_j} \mathbf{Y}^i$ is the empirical mean of $\mathbf{Y}$ on block $B_j$, we can show that Theorem 2.1 can be extended to (regular) *V*-fold cross-validation weights with the following resampling constants:

$$A_W = \frac{2}{V}, \qquad B_W = \frac{1}{\sqrt{V-1}}, \qquad C_W = \frac{\sqrt{n}}{V-1}, \qquad D_W = 1.$$

Additionally, when $V$ does not divide $n$ and the blocks are no longer regular, Theorem 2.1 can also be generalized, but the constants have more complex expressions (see Section 10.7.5 in [1] for details). With *V*-fcv weights, the complexity index is only $V$, but we lose a factor $[(n-1)/(V-1)]^{1/2}$ in the accuracy index. With regard to the accuracy/complexity trade-off, the most accurate cross-validation weights are leave-one-out ($V = n$), whereas the 2-fcv weights are the best from the computational viewpoint (but also the least accurate). The choice of $V$ is thus a trade-off between these two terms and depends on the particular constraints of each problem.

However, it is worth noting that as far as the bound of inequality (7) is concerned, it is not necessarily indispensable to aim for an accuracy index close to 1. Namely, this will result in a corresponding deviation term of order $n^{-1}$, while there is, additionally, another unavoidable deviation term or order $n^{-1/2}$ in the bound. This suggests that an accuracy index of order $o(n^{1/2})$ would actually be sufficient (as $n$ grows large). In other words, using *V*-fcv with $V$ "large" [e.g., $V = \Theta(\log(n))$] would result in only a negligible loss of overall accuracy as compared to leave-one-out. Of course, this discussion is specific to the form of the bound (7). We cannot formally exclude the possibility that a different approach could lead to a different conclusion, unless it can be proven that the deviation terms in (7) cannot be significantly improved, an issue we do not address here.

2.5.2. *Monte Carlo approximation.* When using a Monte Carlo approximation to evaluate $\mathbb{E}_W[\phi(\overline{\mathbf{Y}}^{\langle W - \overline{W} \rangle})]$, we randomly draw a number $B$ of i.i.d. weight vectors $W^1, \ldots, W^B$ and compute $\frac{1}{B} \sum_{j=1}^{B} \phi(\overline{\mathbf{Y}}^{\langle W^j - \overline{W^j} \rangle})$. This method is quite standard in the bootstrap literature and can be improved in several ways (see, e.g., [15], Appendix II).

On the one hand, the number $B$ of draws of $W$ should be taken small enough so that $B$ times the computational cost of evaluating $\phi(\overline{\mathbf{Y}}^{\langle W^j - \overline{W^j} \rangle})$ is still tractable. On the other hand, the number $B$ should be taken large enough to make the Monte Carlo approximation accurate. In our framework, this is quantified more precisely by the following proposition (for bounded weights).

PROPOSITION 2.7. *Let $B \geq 1$ and $W^1, \ldots, W^B$ be i.i.d. exchangeable resampling weight vectors such that $W_1^1 - \overline{W}^1 \in [c_1, c_2]$ a.s. Let $p \in [1, \infty]$ and $\phi : \mathbb{R}^K \to \mathbb{R}$ be any subadditive function bounded by the $\ell_p$-norm. If $\mathbf{Y}$ is a fixed sample, then, for every $\beta \in (0,1)$,*

$$(19) \quad \mathbb{E}_W[\phi(\overline{\mathbf{Y}}^{\langle W - \overline{W} \rangle})] \leq \frac{1}{B} \sum_{j=1}^{B} \phi(\overline{\mathbf{Y}}^{\langle W^j - \overline{W}^j \rangle}) + (c_2 - c_1)\sqrt{\frac{\log(\beta^{-1})}{2B}} \|\widetilde{\sigma}\|_p$$

*holds with probability at least $1 - \beta$, where $\widetilde{\sigma}$ denotes the vector of average absolute deviations to the median, $\widetilde{\sigma} := ((\frac{1}{n}\sum_{i=1}^{n}|\mathbf{Y}_k^i - M_k|))_{1 \leq k \leq K}$ [$M_k$ denoting a median of $(\mathbf{Y}_k^i)_{1 \leq i \leq n}$].*

As a consequence, Proposition 2.7 suggests an explicit correction of the concentration thresholds taking into account $B$ bounded weight vectors. For instance, with Rademacher weights, we can use (19) with $c_2 - c_1 = 2$ and $\beta = \gamma\alpha$ [$\gamma \in (0,1)$]. Then, in the thresholds built from Theorem 2.1 and Proposition 2.2, one can replace $\mathbb{E}_W[\phi(\overline{\mathbf{Y}}^{\langle W - \overline{W} \rangle})]$ by its Monte Carlo approximation at the cost of changing $\alpha$ into $(1-\gamma)\alpha$ and adding $B_W^{-1}\sqrt{\frac{2\log((\gamma\alpha)^{-1})}{B}}\|\widetilde{\sigma}\|_p$ to the threshold.

As $n$ grows large, this remainder term is negligible in comparison to the main one when $B$ is (for instance) of order $n^2$. In practical applications, $B$ can be chosen as a function of $\mathbf{Y}$ because (19) holds conditionally on the observed sample. Therefore, we can use the following strategy: first, compute a rough estimate $t_{\text{est},\alpha}$ of the final threshold [e.g., if $\phi = \|\cdot\|_\infty$ and $\mathbf{Y}$ is Gaussian, take the Bonferroni threshold (10)]. Then choose $B \gg t_{\text{est},\alpha}^2 \|\widetilde{\sigma}\|_p^2 \log((\gamma\alpha)^{-1})$.

## 3. Confidence region using resampled quantiles.

3.1. *Main result.* In this section, we consider a different approach to constructing confidence regions, directly based on the estimation of the quantile via resampling. Once again, since we aim for a nonasymptotic result for $K \gg n$, the standard asymptotic approaches cannot be applied here. For this reason, we base the proposed results on ideas coming from exact randomized tests and consider here the case where $\mathbf{Y}^1$ has a symmetric distribution and where $W$ is an i.i.d. Rademacher weight vector, that is, weights are i.i.d. with $\mathbb{P}(W_i = 1) = \mathbb{P}(W_i = -1) = 1/2$.

The resampling idea applied here is to approximate the quantiles of the distribution $\mathcal{D}(\phi(\overline{\mathbf{Y}} - \mu))$ by the quantiles of the corresponding resampling-based distribution:

$$\mathcal{D}(\phi(\overline{\mathbf{Y}}^{\langle W - \overline{W} \rangle}) | \mathbf{Y}) = \mathcal{D}(\phi((\overline{\mathbf{Y} - \overline{\mathbf{Y}}})^{\langle W \rangle}) | \mathbf{Y}).$$

For this, we take advantage of the symmetry of each $\mathbf{Y}^i$ around its mean. For a function $\phi$, let us define the resampled empirical quantile by

(20) $$q_\alpha(\phi, \mathbf{Y}) := \inf\{x \in \mathbb{R} | \mathbb{P}_W(\phi(\overline{\mathbf{Y}}^{\langle W \rangle}) > x) \leq \alpha\}.$$

The following lemma, close in spirit to exact test results, is easily derived from the "symmetrization trick," that is, from taking advantage of the distribution invariance of the data via sign reversal.

LEMMA 3.1. *Let $\mathbf{Y}$ be a data sample satisfying assumption* (SA) *and $\phi : \mathbb{R}^K \to \mathbb{R}$ be a measurable function. The following then holds:*

(21) $$\mathbb{P}(\phi(\overline{\mathbf{Y}} - \mu) > q_\alpha(\phi, \mathbf{Y} - \mu)) \leq \alpha.$$

Of course, since $q_\alpha(\phi, \mathbf{Y} - \mu)$ still depends on the unknown $\mu$, we cannot use this threshold to get a confidence region of the form (1). It is, in principle, possible to build a confidence region directly from Lemma 3.1 by using the duality between tests and confidence regions, but this would be difficult to compute and not of the desired form (1). Therefore, following the general philosophy of resampling, we propose replacing the true mean $\mu$ by the empirical mean $\overline{\mathbf{Y}}$ in the quantile $q_\alpha(\phi, \mathbf{Y} - \mu)$. The following main technical result of this section gives a nonasymptotic bound on the cost of performing this operation.

THEOREM 3.2. *Fix $\delta, \alpha_0 \in (0, 1)$. Let $\mathbf{Y}$ be a data sample satisfying assumption* (SA). *Let $f : (\mathbb{R}^K)^n \to [0, \infty)$ be a nonnegative function. Let $\phi : \mathbb{R}^K \to \mathbb{R}$ be a nonnegative, subadditive and positive homogeneous function. Define $\widetilde{\phi}(x) := \max(\phi(x), \phi(-x))$. The following holds:*

(22) $$\mathbb{P}(\phi(\overline{\mathbf{Y}} - \mu) > q_{\alpha_0(1-\delta)}(\phi, \mathbf{Y} - \overline{\mathbf{Y}}) + \gamma_1(\alpha_0 \delta) f(\mathbf{Y}))$$
$$\leq \alpha_0 + \mathbb{P}(\widetilde{\phi}(\overline{\mathbf{Y}} - \mu) > f(\mathbf{Y})),$$

*where $\gamma_1(\eta) := \frac{2\overline{\mathcal{B}}(n, \eta/2) - n}{n}$ and $\overline{\mathcal{B}}(n, \eta) := \max\{k \in \{0, \ldots, n\} | 2^{-n} \sum_{i=k}^n \binom{n}{i} \geq \eta\}$ is the upper quantile function of a* Binomial$(n, \frac{1}{2})$ *variable.*

In this result, the resampled quantile term $q_{\alpha_0(1-\delta)}(\phi, \mathbf{Y} - \overline{\mathbf{Y}})$ should be interpreted as the main term of the threshold and the rest, involving the function $f$, as a remainder term. In the usual resampling philosophy, one would only consider the main term at the target level, that is, $\alpha_0 = \alpha$ and $\delta =$

0. Here, the additional remainder terms are introduced to account rigorously for the validity of the result in a nonasymptotic setting. These remainder terms have two effects: first, the resampled quantile in the main term is computed at a "shrunk" error level $\alpha_0(1-\delta) < \alpha$ and, secondly, there is an additional additive term in the threshold itself.

The role of the parameters $\delta$, $\alpha_0$ and $f$ is to strike a balance between these effects. Generally speaking, $f$ should be an available upper bound on a quantile of $\widetilde{\phi}(\overline{\mathbf{Y}} - \mu)$ at a level $\alpha_1 \ll \alpha_0$. On the left-hand side, $f$ appears in the threshold with the factor $\gamma_1$, which can be more explicitly bounded by

$$(23) \qquad \gamma_1(\alpha_0 \delta) \leq \left(\frac{2\log(2/(\alpha_0 \delta))}{n}\right)^{1/2},$$

using Hoeffding's inequality. The above result therefore transforms a possibly coarse "a priori" bound $f$ on quantiles into a more accurate quantile bound based on a main term estimated by resampling and a remainder term based on $f$ multiplied by a small factor.

In order to get a clearer insight, let us consider an example of specific choices for the parameters $\delta, \alpha_0$ and $f$ in the Gaussian case. First, choose $\delta = \Theta(n^{-\gamma})$ and $\frac{\alpha_0}{\alpha} = 1 - \Theta(n^{-\gamma})$ for some $\gamma > 0$, say $\gamma = 1$. This way, the main term is the resampled quantile at level $\alpha_0(1-\delta) = \alpha(1 - \Theta(n^{-\gamma}))$. For the choice of $f$, let us choose Bonferroni's threshold (10) at level $\alpha_1 = (\alpha - \alpha_0) = \Theta(n^{-\gamma})$ so that the overall probability control in (22) is really at the target level $\alpha$. Then $f_{\mathrm{Bonf}}(\mathbf{Y}) \leq \Theta((\log(Kn^\gamma)/n)^{1/2})$ and, using (23), we conclude that the remainder term is bounded by $\Theta(\log(Kn^\gamma)/n)$. This is indeed a remainder term with respect to the main term which is of order at least $\Theta(n^{-1/2})$ as $n$ grows [assuming that the dimension $K(n)$ grows subexponentially with $n$].

There are other possibilities for choosing $f$, depending on the context: the Bonferroni threshold can be correspondingly adapted to the non-Gaussian case when an upper bound on the tail of each coordinate is available. This still makes the remainder term directly dependent on $K$ and a possibly more interesting idea is to recycle the results of Section 2 (when the data is either Gaussian or bounded and symmetric) and plug in the thresholds derived there for the function $f$.

Finally, if the a priori bound on the quantiles is too coarse, it is possible to iterate the process and estimate smaller quantiles more accurately by again using resampling. Namely, by iteration of Theorem 3.2, we obtain the following corollary.

COROLLARY 3.3. *Fix $J$ a positive integer, $(\alpha_i)_{i=0,\ldots,J-1}$ a finite sequence in $(0,1)$ and $\delta \in (0,1)$. Consider $\mathbf{Y}$, $f$, $\phi$ and $\widetilde{\phi}$ as in Theorem 3.2. The*

*following then holds:*

$$
\begin{aligned}
\mathbb{P}\Bigg(\phi(\overline{\mathbf{Y}} - \mu) &> q_{\alpha_0(1-\delta)}(\phi, \mathbf{Y} - \overline{\mathbf{Y}}) \\
&+ \sum_{i=1}^{J-1} \gamma_i q_{\alpha_i(1-\delta)}(\widetilde{\phi}, \mathbf{Y} - \overline{\mathbf{Y}}) + \gamma_J f(\mathbf{Y})\Bigg) \\
&\leq \sum_{i=0}^{J-1} \alpha_i + \mathbb{P}(\widetilde{\phi}(\overline{\mathbf{Y}} - \mu) > f(\mathbf{Y})),
\end{aligned}
$$
(24)

*where, for* $k \geq 1$, $\gamma_k := n^{-k} \prod_{i=0}^{k-1}(2\overline{\mathcal{B}}(n, \frac{\alpha_i \delta}{2}) - n)$.

The rationale behind this result is that the sum appearing inside the probability in (24) should be interpreted as a series of corrective terms of decreasing order of magnitude because we expect the sequence $\gamma_k$ to be sharply decreasing. From (23), this will be the case if the levels are such that $\alpha_i \gg \exp(-n)$.

The conclusion is that even if the a priori available bound $f$ on small quantiles is not sharp, its contribution to the threshold can be made small in comparison to the (more accurate) resampling terms. The counterpart to be paid is the loss in the level and the additional terms in the threshold; for large $n$, these terms decay very rapidly, but for small $n$, they may still result in a nonnegligible contribution; in this case, a precise tuning of the parameters $J, (\alpha_i), \delta$ and $f$ is of much more importance and also more delicate.

At this point, we should also mention that the remainder terms given by Theorem 3.2 and Corollary 3.3 are certainly overestimated, even if $f$ is very well chosen. This makes the theoretical thresholds slightly too conservative in general (particularly for small values of $n$). From simulations not reported here (see [2] and Section 4.3 below), it even appears that the remainder terms could be (almost) unnecessary in standard situations, even for $n$ relatively small. Proving this fact rigorously in a nonasymptotic setting, possibly with some additional assumption on the distribution of $\mathbf{Y}$, remains an open issue. Another interesting open problem would be to obtain a self-contained result based on the symmetry assumption (SA) alone [or a negative result proving that (SA) is not sufficient for a distribution-free result of this form].

3.2. *Practical computation of the resampled quantile.* Since the above results use Rademacher weight vectors, the exact computation of the quantile $q_\alpha$ requires, in principle, $2^n$ iterations and is thus too complex as $n$ becomes large. Parallel to what was proposed for the concentration-based thresholds in Section 2.5, one can, as a first solution, consider a blockwise Rademacher resampling scheme or, equivalently, applying the previous method to a block-averaged sample, at the cost of a (possibly substantial) loss in accuracy.

A possibly better way to address this issue is by means of Monte Carlo quantile approximation, on which we now focus. Let $\mathbf{W}$ denote an $n \times B$ matrix of i.i.d. Rademacher weights (independent of all other variables) and define

$$\widetilde{q}_\alpha(\phi, \mathbf{Y}, \mathbf{W}) := \inf\left\{x \in \mathbb{R} \,\bigg|\, \frac{1}{B}\sum_{j=1}^{B} \mathbb{1}\{\phi(\overline{\mathbf{Y}}^{\langle \mathbf{W}^j \rangle}) \geq x\} \leq \alpha\right\},$$

that is, $\widetilde{q}_\alpha$ is defined in the same way as $q_\alpha$, except that the true distribution $\mathbb{P}_W$ of the Rademacher weight vector is replaced by the empirical distribution constructed from the columns of $\mathbf{W}$, $\widetilde{\mathbb{P}}_\mathbf{W} = B^{-1}\sum_{j=1}^{B}\delta_{\mathbf{W}^j}$; note that the strict inequality $\phi(\overline{\mathbf{Y}}^{\langle \mathbf{W} \rangle}) > x$ in (20) was replaced by $\phi(\overline{\mathbf{Y}}^{\langle \mathbf{W}^j \rangle}) \geq x$ for technical reasons. The following result then holds.

PROPOSITION 3.4. *Consider the same conditions as in Theorem 3.2, except that the function $f$ can now be a function of both $\mathbf{Y}$ and $\mathbf{W}$. We have*

$$\mathbb{P}_{\mathbf{Y},\mathbf{W}}(\phi(\overline{\mathbf{Y}} - \mu) > \widetilde{q}_{\alpha_0(1-\delta)}(\phi, \mathbf{Y} - \overline{\mathbf{Y}}, \mathbf{W}) + \gamma(\mathbf{W}, \alpha_0\delta)f(\mathbf{Y}, \mathbf{W}))$$
$$\leq \widetilde{\alpha}_0 + \mathbb{P}_{\mathbf{Y},\mathbf{W}}(\widetilde{\phi}(\overline{\mathbf{Y}} - \mu) > f(\mathbf{Y}, \mathbf{W})),$$

*where $\widetilde{\alpha}_0 := \frac{\lfloor B\alpha_0 \rfloor + 1}{B+1} \leq \alpha_0 + \frac{1}{B+1}$ and $\gamma(\mathbf{W}, \eta) := \max\{y \geq 0 | \frac{1}{B}\sum_{j=1}^{B} \mathbb{1}\{|\overline{\mathbf{W}}^j| \geq y\} \geq \eta\}$ is the $(1-\eta)$-quantile of $|\overline{W}|$ under the empirical distribution $\widetilde{\mathbb{P}}_\mathbf{W}$.*

Note that, for practical purposes, we can choose $f(\mathbf{W}, \mathbf{Y})$ to depend on $\mathbf{Y}$ only and use another type of bound to control the last term on the right-hand side, as in the earlier discussion. The above result tells us that if, in Theorem 3.2, we replace the true quantile by an empirical quantile based on $B$ i.i.d. weight vectors and the factor $\gamma_1$ is similarly replaced by an empirical quantile of $|\overline{W}|$, then we lose at most $(B+1)^{-1}$ in the corresponding covering probability. Furthermore, it can easily be seen that if $\alpha_0$ is taken to be a positive multiple of $(B+1)^{-1}$, then there is no loss in the final covering probability (i.e., $\widetilde{\alpha}_0 = \alpha_0$).

## 4. Discussion and concluding remarks.

4.1. *Estimating $\|\sigma\|_p$*. In the concentration approach and in the Gaussian case, the derived thresholds depend explicitly on the $\ell_p$-norm of the vector of standard deviations $\sigma = (\sigma_k)_k$ (an upper bound on this quantity can also be used). While we have left aside the problem of determining this parameter if no prior information is available, it is possible to estimate $\sigma$ by its empirical counterpart

$$\widehat{\sigma} := \left(\sqrt{\frac{1}{n}\sum_{i=1}^{n}(\mathbf{Y}_k^i - \overline{\mathbf{Y}}_k)^2}\right)_{1 \leq k \leq K}.$$

Interestingly, the quantity $\|\widehat{\sigma}\|_p$ enjoys the same type of concentration property as the resampled expectations considered in Section 2.3 so that we can derive, by a similar argument, a *dimension-free* confidence bound for $\|\sigma\|_p$, as follows.

PROPOSITION 4.1. *Assume that* $\mathbf{Y}$ *satisfies* (GA). *Then, with probability at least* $1 - \delta$,

$$\|\sigma\|_p \leq \left(C_n - \frac{1}{\sqrt{n}}\overline{\Phi}^{-1}\left(\frac{\delta}{2}\right)\right)^{-1}\|\widehat{\sigma}\|_p, \tag{25}$$

*where* $C_n = \sqrt{\frac{2}{n}}\frac{\Gamma(n/2)}{\Gamma((n-1)/2)}$.

It can easily be checked via Stirling's formula that $C_n = 1 - O(n^{-1})$, so replacing $\|\sigma\|_p$ by the above upper bound does not make the corresponding thresholds significantly more conservative.

A similar question holds for the parameter $M$ in the bounded case. In practical applications, an absolute bound on the possible data values is often known (e.g., from physical or biological constraints). It can also be estimated, but it seems harder to obtain a rigorous nonasymptotic control on the level of the resulting threshold in the general bounded case.

A different, and potentially more important, problem arises if the vector of variances $\sigma$ is not constant. Since the confidence regions proposed in this paper are isotropic, they will—inevitably—tend to be conservative when the variances of the coordinates are very different. The standard way to address this issue is to consider studentized data. While this would solve this heteroscedasticity issue, it also renders void the assumption of independent data points—a crucial assumption in all of our proofs. Therefore, generalizing our results to studentized observations is an important, but probably challenging, direction for future work.

4.2. *Interpretation and use of $\phi$-confidence regions.* We have built high-dimensional confidence regions taking the form of "$\phi$-balls" [where $\phi$ can be any $\ell^p$-norm with $p \geq 1$, but more general choices are possible, such as $\phi(x) = \sup_k (x_k)_+$]. Such confidence regions in very high dimension are certainly quite difficult to visualize and one can ask how they are to be interpreted. In our opinion, the most intuitive and interesting interpretation again comes from learning theory, by regarding $\phi$ as a type of loss function. In this sense, a $\phi$-confidence region is an upper confidence bound on some relevant loss measure of the estimator $\mathbf{Y}$ to the target $\mu$. Additionally, in the particular case when $\phi = \sup_k (x_k)_+$ or $\phi = \|\cdot\|_\infty$, the corresponding regions can be interpreted as simultaneous confidence intervals over all coordinate means.

The results presented here can also provide confidence intervals for the $\ell^p$-risk (i.e., the averaged $\phi$-loss) for the estimator $\overline{\mathbf{Y}}$ of the mean vector $\mu$. Indeed, combining (12) and Proposition 2.5(ii), we derive that for a Gaussian sample $\mathbf{Y}$ and any $p \in [1, \infty]$, the upper bound

$$(26) \qquad \mathbb{E}\|\overline{\mathbf{Y}} - \mu\|_p < \frac{\mathbb{E}_W[\|\overline{\mathbf{Y}}^{\langle W - \overline{W}\rangle}\|_p]}{B_W} + \frac{\|\sigma\|_p C_W}{nB_W} \overline{\Phi}^{-1}(\alpha/2),$$

holds with probability at least $1 - \alpha$ and a similar lower bound holds. It is worth noting that the rate $C_W/(nB_W)$ is close to $n^{-1}$ for most of the weights, meaning that resampling provides a much better estimate of $\mathbb{E}\|\overline{\mathbf{Y}} - \mu\|_p$ than $\|\overline{\mathbf{Y}} - \mu\|_p$ itself. This stabilization effect of resampling is well known in standard asymptotic settings (see, e.g., [15]).

The $\ell^p$-risk is also related to the leave-one-out estimation of the prediction risk. Indeed, consider using $\overline{\mathbf{Y}}$ for *predicting* a new data point $\mathbf{Y}^{n+1} \sim \mathbf{Y}^1$ [independent of $\mathbf{Y} = (\mathbf{Y}^1, \ldots, \mathbf{Y}^n)$]. The corresponding $\ell^p$-prediction risk is given by $\mathbb{E}\|\overline{\mathbf{Y}} - \mathbf{Y}^{n+1}\|_p$. In the Gaussian setting, this prediction risk is proportional to the $\ell^p$-risk: $\mathbb{E}\|\overline{\mathbf{Y}} - \mu\|_p = (n+1)^{1/2}\mathbb{E}\|\overline{\mathbf{Y}} - \mathbf{Y}^{n+1}\|_p$, so the previous resampling estimator of the $\ell^p$-risk also leads to an estimator of the prediction risk. In particular, using leave-one-out weights and denoting by $\overline{\mathbf{Y}}^{(-i)}$, the mean of the $(\mathbf{Y}^j, j \neq i, 1 \leq j \leq n)$, our results prove that the leave-one-out estimator

$$\frac{1}{n} \sum_{i=1}^n \|\overline{\mathbf{Y}}^{(-i)} - \mathbf{Y}^i\|_p$$

correctly estimates the prediction risk [up to the factor $(1 - 1/n^2)^{1/2} \simeq 1$].

Finally, another important field of application is hypothesis testing. When $\phi = \sup_k (x_k)_+$ or $\phi = \|\cdot\|_\infty$, the thresholds derived here can be used to derive *multiple testing* procedures for the value of the mean of each coordinate. This question is extensively developed in the companion paper [2]. It is also possible to take advantage of the generality of our results, where $\phi$ is allowed to be any $\ell^p$-norm with $p \geq 1$, for single global hypothesis testing. The confidence regions can be used straightforwardly to test several single global hypotheses, such as $\mu = \mu^\star$ against $\|\mu - \mu^\star\|_p \geq R > 0$. Depending on particular features of the problem, having the choice between different functions $\phi$ allows us to take into account specific forms of alternative hypotheses in the construction of the threshold.

4.3. *Simulation study.* In the companion paper [2] (Section 4), a simulation study compares the thresholds built in this paper and Bonferroni's threshold, using $\phi = \|\cdot\|_\infty$, considering Gaussian data with different levels of correlations and assuming the coordinate variance $\sigma$ to be constant and known. Without entering into details, its general conclusions are as follows. First, all of the thresholds proposed in the present paper can improve on Bonferroni's when the correlations are strong enough. Even though our thresholds are seen to be more conservative than the "ideal" one (i.e., the true quantile), they all exhibit adaptivity to the correlations, as expected from their construction. However, when the vector coordinates are close to being independent, the proposed thresholds are somewhat more conservative than Bonferroni's (the latter being essentially optimal in this case).

The second observation made on the simulations is that the quantile approach generally appears to be less conservative than the concentration approach. However, the remaining advantage of the concentration approach is that it can be combined with Bonferroni's threshold (using Proposition 2.2) so that one can almost take "the less conservative of the two" and only suffer a negligible loss if the Bonferroni threshold turns out to be better. Also, recall that the concentration threshold can be of use for the remainder terms of the quantile threshold.

Finally, we also tested the resampled quantile without remainder term (i.e., taking the raw resampled quantile of the empirically centered data at the desired level, without modification). Although this threshold is not theoretically justified in the present work, it appeared to be very close to the ideal threshold in the performed simulations. This supports the conjecture that the remainder terms in the theoretical threshold could either be made significantly smaller or, possibly, even completely dropped in some cases.

4.4. *Comparing nonasymptotic and asymptotic approaches.* Although simulations have shown that the various thresholds proposed here can outperform Bonferroni's when significant correlations are present, we have also noticed that these thresholds are generally noticeably more conservative than the ideal ones (the true quantiles), especially for small values of $n$. Moreover, taking into account other sources of error such as the estimation of $\|\sigma\|_p$ as above, or Monte Carlo approximations, will result in even more conservative thresholds. The main reason for this additional conservativeness is that our control on the level is *nonasymptotic*, that is, valid for every fixed $K$ and $n$. In this sense, it would be somewhat unfair to compare the thresholds proposed here to those of "traditional" resampling theory that are only proved to be valid asymptotically in $n$ and for fixed $K$. The nonasymptotic results derived here can nevertheless also be used for an asymptotic analysis, in a setting where $K(n)$ is a function of $n$, and possibly rapidly (say, exponentially) growing. This type of situation seems to have been only scarcely touched by existing asymptotic approaches. In this sense, in practical situations, we can envision "cheating" somewhat and replacing the theoretical thresholds by their leading component [under some mild assumptions on the growth of $K(n)$] as $n$ tends to infinity. From a theoretical point of view, an interesting avenue for future endeavors is to prove that the thresholds considered here, while certainly not second order correct, are at least asymptotically optimal under various dependency conditions.

## 5. Proofs.

5.1. *Confidence regions using concentration.* In this section, we prove all of the statements of Section 2 except computations of resampling weight constants (made in Section 5.3).

5.1.1. *Comparison in expectation.*

PROOF OF PROPOSITION 2.3. Denoting by $\boldsymbol{\Sigma}$ the common covariance matrix of the $\mathbf{Y}^i$, we have $\mathcal{D}(\overline{\mathbf{Y}}^{\langle W-\overline{W}\rangle}|W) = \mathcal{N}(0, (n^{-1}\sum_{i=1}^n (W_i-\overline{W})^2)n^{-1}\boldsymbol{\Sigma})$ and the result follows because $\mathcal{D}(\overline{\mathbf{Y}} - \mu) = \mathcal{N}(0, n^{-1}\boldsymbol{\Sigma})$ and $\phi$ is positive homogeneous. This last assumption is, of course, unnecessary if it holds that $\sum_{i=1}^n (W_i - \overline{W})^2 = n$ a.s. □

PROOF OF PROPOSITION 2.4. By independence between $W$ and $\mathbf{Y}$, exchangeability of $W$ and the positive homogeneity of $\phi$, for every realization of $\mathbf{Y}$, we have

$$A_W \phi(\overline{\mathbf{Y}} - \mu) = \phi\left(\mathbb{E}\left[\frac{1}{n}\sum_{i=1}^n |W_i - \overline{W}|(\mathbf{Y}^i - \mu)\Big|\mathbf{Y}\right]\right).$$

Then, by convexity of $\phi$,

$$A_W \phi(\overline{\mathbf{Y}} - \mu) \leq \mathbb{E}\left[\phi\left(\frac{1}{n}\sum_{i=1}^n |W_i - \overline{W}|(\mathbf{Y}^i - \mu)\right)\Big|\mathbf{Y}\right].$$

We integrate with respect to $\mathbf{Y}$ and use the symmetry of the $\mathbf{Y}^i$ with respect to $\mu$ and, again, the independence between $W$ and $\mathbf{Y}$ to show, finally, that

$$A_W \mathbb{E}[\phi(\overline{\mathbf{Y}} - \mu)] \leq \mathbb{E}\left[\phi\left(\frac{1}{n}\sum_{i=1}^n |W_i - \overline{W}|(\mathbf{Y}^i - \mu)\right)\right]$$

$$= \mathbb{E}\left[\phi\left(\frac{1}{n}\sum_{i=1}^n (W_i - \overline{W})(\mathbf{Y}^i - \mu)\right)\right] = \mathbb{E}[\phi(\overline{\mathbf{Y}}^{\langle W-\overline{W}\rangle})].$$

The point (ii) is proved via the following chain of inequalities:

$$\mathbb{E}[\phi(\overline{\mathbf{Y}}^{\langle W-\overline{W}\rangle})] \leq \mathbb{E}\left[\phi\left(\frac{1}{n}\sum_{i=1}^{n}(W_i - x_0)(\mathbf{Y}^i - \mu)\right)\right]$$
$$+ \mathbb{E}\left[\phi\left(\frac{1}{n}\sum_{i=1}^{n}(x_0 - \overline{W})(\mathbf{Y}^i - \mu)\right)\right]$$
$$= \mathbb{E}\left[\phi\left(\frac{1}{n}\sum_{i=1}^{n}|W_i - x_0|(\mathbf{Y}^i - \mu)\right)\right]$$
$$+ \mathbb{E}\left[\phi\left(\frac{1}{n}\sum_{i=1}^{n}|x_0 - \overline{W}|(\mathbf{Y}^i - \mu)\right)\right]$$
$$\leq (a + \mathbb{E}|\overline{W} - x_0|)\mathbb{E}[\phi(\overline{\mathbf{Y}} - \mu)].$$

In the second line, we used, as before, the symmetry of the $\mathbf{Y}^i$ with respect to $\mu$, together with the independence of $W$ and $\mathbf{Y}$. In the last inequality, we used the assumption $|W_i - x_0| = a$ a.s. and the positive homogeneity of $\phi$. □

### 5.1.2. *Concentration inequalities.*

PROOF OF PROPOSITION 2.5. Here, we use concentration principles following closely the approach in [24], Section 3.2.4. The essential ingredient is the Gaussian concentration theorem of Cirel'son, Ibragimov and Sudakov ([7] and recalled in [24], Theorem 3.8), stating that if $F$ is a Lipschitz function on $\mathbb{R}^N$ with constant $L$, then, for the standard Gaussian measure on $\mathbb{R}^N$, we have $\mathbb{P}(F \geq \mathbb{E}[F] + t) \leq 2\overline{\Phi}(t/L)$.

Let us denote by $\mathbf{A}$ a square root of the common covariance matrix of the $\mathbf{Y}^i$. If $\zeta_i$ is a $K$-dimensional, standard normal vector, then $\mathbf{A}\zeta_i$ has the same distribution as $\mathbf{Y}^i - \mu$. For all $\zeta \in (\mathbb{R}^K)^n$, we let $T_1(\zeta) := \phi(\frac{1}{n}\sum_{i=1}^{n}\mathbf{A}\zeta_i)$ and $T_2(\zeta) := \mathbb{E}[\phi(\frac{1}{n}\sum_{i=1}^{n}(W_i - \overline{W})\mathbf{A}\zeta_i)]$. If we endow $(\mathbb{R}^K)^n$ with the standard Gaussian measure, then $T_1$ (resp., $T_2$) has the same distribution as $\phi(\overline{\mathbf{Y}} - \mu)$ [resp., $\phi(\overline{\mathbf{Y}}^{\langle W-\overline{W}\rangle})$].

From the Gaussian concentration theorem recalled above, in order to reach the conclusion, we therefore only need to establish that $T_1$ (resp., $T_2$) is a Lipschitz function with constant $\|\sigma\|_p/\sqrt{n}$ (resp., $\|\sigma\|_p C_W/n$) with respect to the Euclidean norm $\|\cdot\|_{2,Kn}$ on $(\mathbb{R}^K)^n$. Let $\zeta, \zeta' \in (\mathbb{R}^K)^n$ and denote by $(a_k)_{1\leq k\leq K}$ the rows of $\mathbf{A}$. Using the fact that $\phi$ is 1-Lipschitz with respect to the $\ell_p$-norm (because it is subadditive and bounded by the $\ell_p$-norm), we

get

$$|T_1(\zeta) - T_1(\zeta')| \leq \left\|\frac{1}{n}\sum_{i=1}^n \mathbf{A}(\zeta_i - \zeta'_i)\right\|_p = \left\|\left(\left\langle a_k, \frac{1}{n}\sum_{i=1}^n (\zeta_i - \zeta'_i)\right\rangle\right)_k\right\|_p.$$

For each coordinate $k$, by the Cauchy–Schwarz inequality and since $\|a_k\|_2 = \sigma_k$, we deduce that

$$\left|\left\langle a_k, \frac{1}{n}\sum_{i=1}^n (\zeta_i - \zeta'_i)\right\rangle\right| \leq \sigma_k \left\|\frac{1}{n}\sum_{i=1}^n (\zeta_i - \zeta'_i)\right\|_2.$$

Therefore, we get

$$|T_1(\zeta) - T_1(\zeta')| \leq \|\sigma\|_p \left\|\frac{1}{n}\sum_{i=1}^n (\zeta_i - \zeta'_i)\right\|_2 \leq \frac{\|\sigma\|_p}{\sqrt{n}}\|\zeta - \zeta'\|_{2,Kn},$$

using the convexity of $x \in \mathbb{R}^K \mapsto \|x\|_2^2$, and we obtain (i). For $T_2$, we use the same method as for $T_1$ to obtain

(27)
$$|T_2(\zeta) - T_2(\zeta')| \leq \|\sigma\|_p \mathbb{E}\left\|\frac{1}{n}\sum_{i=1}^n (W_i - \overline{W})(\zeta_i - \zeta'_i)\right\|_2$$

$$\leq \frac{\|\sigma\|_p}{n}\sqrt{\mathbb{E}\left\|\sum_{i=1}^n (W_i - \overline{W})(\zeta_i - \zeta'_i)\right\|_2^2}.$$

Note that since $(\sum_{i=1}^n (W_i - \overline{W}))^2 = 0$, we have $\mathbb{E}(W_1 - \overline{W})(W_2 - \overline{W}) = -C_W^2/n$. We now develop $\|\sum_{i=1}^n (W_i - \overline{W})(\zeta_i - \zeta'_i)\|_2^2$ in the Euclidean space $\mathbb{R}^K$:

$$\mathbb{E}\left\|\sum_{i=1}^n (W_i - \overline{W})(\zeta_i - \zeta'_i)\right\|_2^2$$

$$= C_W^2(1 - n^{-1})\sum_{i=1}^n \|\zeta_i - \zeta'_i\|_2^2 - \frac{C_W^2}{n}\sum_{i\neq j}\langle \zeta_i - \zeta'_i, \zeta_j - \zeta'_j\rangle$$

$$= C_W^2 \sum_{i=1}^n \|\zeta_i - \zeta'_i\|_2^2 - \frac{C_W^2}{n}\left\|\sum_{i=1}^n (\zeta_i - \zeta'_i)\right\|_2^2.$$

Consequently,

(28) $$\mathbb{E}\left\|\sum_{i=1}^n (W_i - \overline{W})(\zeta_i - \zeta'_i)\right\|_2^2 \leq C_W^2 \sum_{i=1}^n \|\zeta_i - \zeta'_i\|_2^2 = C_W^2 \|\zeta - \zeta'\|_{2,Kn}^2.$$

Combining expression (27) and (28), we find that $T_2$ is $\|\sigma\|_p C_W/n$-Lipschitz.
□

REMARK 5.1. The proof of Proposition 2.5 is still valid under the weaker assumption (instead of exchangeability of $W$) that $\mathbb{E}[(W_i - \overline{W})(W_j - \overline{W})]$ can only take two possible values, depending on whether or not $i = j$.

### 5.1.3. *Main results.*

PROOF OF THEOREM 2.1. The case (BA) $(p, M)$ and (SA) is obtained by combining Propositions 2.4 and 2.6. The (GA) case is a straightforward consequence of Proposition 2.3 and the proof of Proposition 2.5 (considering the Lipschitz function $T_1 - T_2$). □

PROOF OF PROPOSITION 2.2. From Proposition 2.5(i), with probability at least $1 - \alpha(1-\delta)$, $\phi(\overline{\mathbf{Y}} - \mu)$ is less than or equal to the minimum of $t_{\alpha(1-\delta)}$ and $\mathbb{E}[\phi(\overline{\mathbf{Y}} - \mu)] + \frac{\|\sigma\|_p \overline{\Phi}^{-1}(\alpha(1-\delta)/2)}{\sqrt{n}}$ (since both of these thresholds are deterministic). In addition, Propositions 2.3 and 2.5(ii) give that with probability at least $1 - \alpha\delta$, $\mathbb{E}[\phi(\overline{\mathbf{Y}} - \mu)] \leq \frac{\mathbb{E}_W[\phi(\overline{\mathbf{Y}}^{\langle W - \overline{W}\rangle})]}{B_W} + \frac{\|\sigma\|_p C_W}{B_W n} \overline{\Phi}^{-1}(\alpha\delta/2)$. The result follows by combining the last two expressions. □

### 5.1.4. *Monte Carlo approximation.*

PROOF OF PROPOSITION 2.7. The idea of the proof is to apply McDiarmid's inequality (see [25]) conditionally on $\mathbf{Y}$. For any realizations $W$ and $W'$ of the resampling weight vector and any $\nu \in \mathbb{R}^k$, we have

$$|\phi(\overline{\mathbf{Y}}^{\langle W - \overline{W}\rangle}) - \phi(\overline{\mathbf{Y}}^{\langle W' - \overline{W'}\rangle})| \leq \phi(\overline{\mathbf{Y}}^{\langle W - \overline{W}\rangle} - \overline{\mathbf{Y}}^{\langle W' - \overline{W'}\rangle})$$

$$\leq \frac{c_2 - c_1}{n} \left\| \left( \sum_{i=1}^n |\mathbf{Y}_k^i - \nu_k| \right)_k \right\|_p$$

since $\phi$ is subadditive, bounded by the $\ell_p$-norm and $W_i - \overline{W} \in [c_1, c_2]$ a.s.

The sample $\mathbf{Y}$ being deterministic, we can take $\nu_k$ equal to a median $M_k$ of $(\mathbf{Y}_k^i)_{1 \leq i \leq n}$. Since $W^1, \ldots, W^B$ are independent, McDiarmid's inequality gives (19). □

### 5.1.5. *Estimation of the variance.*

PROOF OF PROPOSITION 4.1. We use the same notation and approach based on Gaussian concentration as in the proof of Proposition 2.5. Writing $\mathbf{Y}^i - \mu = \mathbf{A}\zeta_i$, we upper bound the Lipschitz constant of $\|\widehat{\sigma}\|_p$ as a function of $\zeta = (\zeta_1, \ldots, \zeta_n)$: given $\zeta, \zeta' \in (\mathbb{R}^K)^n$, we have

$$\|\widehat{\sigma}(\zeta)\|_p - \|\widehat{\sigma}(\zeta')\|_p \leq \|\widehat{\sigma}(\zeta) - \widehat{\sigma}(\zeta')\|_p$$

$$\leq \left\| \left( \frac{1}{n} \sum_{i=1}^{n} \langle a_k, (\zeta_i - \overline{\zeta}) - (\zeta'_i - \overline{\zeta}') \rangle^2 \right)^{1/2}_k \right\|_p$$

$$\leq \frac{\|\sigma\|_p}{\sqrt{n}} \left( \sum_{i=1}^{n} \|(\zeta_i - \overline{\zeta}) - (\zeta'_i - \overline{\zeta}')\|_2^2 \right)^{1/2}.$$

We then additionally have

$$\sum_{i=1}^{n} \|(\zeta_i - \overline{\zeta}) - (\zeta'_i - \overline{\zeta}')\|_2^2 = \sum_{i=1}^{n} \|\zeta_i - \zeta'_i\|_2^2 - n\|\overline{\zeta} - \overline{\zeta}'\|_2^2 \leq \|\zeta - \zeta'\|_{2,Kn}^2,$$

allowing us to conclude that $\|\widehat{\sigma}(\zeta)\|_p$ has Lipschitz constant $\frac{\|\sigma\|_p}{\sqrt{n}}$. Concerning the expectation, observe that for each coordinate $k$, the variable $\sqrt{n}\widehat{\sigma}_k/\sigma_k$ has the same distribution as the square root of a $\chi^2(n-1)$ variable. Elementary calculations for the expectation of such a variable lead to $\mathbb{E}[\widehat{\sigma}_k] = C_n\sigma_k$. We finally conclude that with probability at least $1-\delta$, the following inequality holds:

$$C_n\|\sigma\|_p = \|\mathbb{E}[\widehat{\sigma}]\|_p \leq \mathbb{E}[\|\widehat{\sigma}\|_p] \leq \|\widehat{\sigma}\|_p + \frac{\|\sigma\|_p}{\sqrt{n}}\overline{\Phi}^{-1}\left(\frac{\delta}{2}\right).$$

Solving this inequality in $\|\sigma\|_p$ yields the result. □

5.2. *Quantiles.* Recall the following inequality coming from the definition of the quantile $q_\alpha$: for any fixed $\mathbf{Y}$,

(29) $\quad \mathbb{P}_W(\phi(\overline{\mathbf{Y}}^{\langle W \rangle}) > q_\alpha(\phi, \mathbf{Y})) \leq \alpha \leq \mathbb{P}_W(\phi(\overline{\mathbf{Y}}^{\langle W \rangle}) \geq q_\alpha(\phi, \mathbf{Y})).$

PROOF OF LEMMA 3.1. We introduce the notation $\mathbf{Y} \bullet W = \mathbf{Y} \cdot \text{diag}(W)$ for the matrix obtained by multiplying the $i$th column of $\mathbf{Y}$ by $W_i$, $i = 1, \ldots, n$. We then have

$$\mathbb{P}_\mathbf{Y}(\phi(\overline{\mathbf{Y}} - \mu) > q_\alpha(\phi, \mathbf{Y} - \mu))$$
(30)
$$= \mathbb{E}_W[\mathbb{P}_\mathbf{Y}(\phi((\overline{\mathbf{Y} - \mu})^{\langle W \rangle}) > q_\alpha(\phi, (\mathbf{Y} - \mu) \bullet W))]$$
$$= \mathbb{E}_\mathbf{Y}[\mathbb{P}_W(\phi((\overline{\mathbf{Y} - \mu})^{\langle W \rangle}) > q_\alpha(\phi, \mathbf{Y} - \mu))] \leq \alpha.$$

The first equality is due to the fact that the distribution of $\mathbf{Y}$ satisfies assumption (SA), hence the distribution of $(\mathbf{Y} - \mu)$ is invariant under multiplying by (arbitrary) signs $W \in \{-1, 1\}^n$. In the second equality, we used Fubini's theorem and the fact that for any arbitrary signs $W$, as above, $q_\alpha(\phi, (\mathbf{Y} - \mu) \bullet W) = q_\alpha(\phi, \mathbf{Y} - \mu)$. Finally, the last inequality follows from (29). □

PROOF OF THEOREM 3.2. Write $\gamma_1 = \gamma_1(\alpha_0 \delta)$ for short and define the event

$$\mathcal{E} := \{\mathbf{Y} | q_{\alpha_0}(\phi, \mathbf{Y} - \mu) \leq q_{\alpha_0(1-\delta)}(\phi, \mathbf{Y} - \overline{\mathbf{Y}}) + \gamma_1 f(\mathbf{Y})\}.$$

We then have, using (30),

$$\begin{aligned}
\mathbb{P}(\phi(\overline{\mathbf{Y}} - \mu) &> q_{\alpha_0(1-\delta)}(\phi, \mathbf{Y} - \overline{\mathbf{Y}}) + \gamma_1 f(\mathbf{Y})) \\
&\leq \mathbb{P}(\phi(\overline{\mathbf{Y}} - \mu) > q_{\alpha_0}(\phi, \mathbf{Y} - \mu)) + \mathbb{P}(\mathbf{Y} \in \mathcal{E}^c) \\
&\leq \alpha_0 + \mathbb{P}(\mathbf{Y} \in \mathcal{E}^c).
\end{aligned} \tag{31}$$

We now concentrate on the event $\mathcal{E}^c$. Using the subadditivity of $\phi$ and the fact that $(\overline{\mathbf{Y} - \mu})^{\langle W \rangle} = (\overline{\mathbf{Y} - \overline{\mathbf{Y}}})^{\langle W \rangle} + \overline{W}(\overline{\mathbf{Y}} - \mu)$, we have, for any fixed $\mathbf{Y} \in \mathcal{E}^c$,

$$\begin{aligned}
\alpha_0 &\leq \mathbb{P}_W(\phi((\overline{\mathbf{Y} - \mu})^{\langle W \rangle}) \geq q_{\alpha_0}(\phi, \mathbf{Y} - \mu)) \\
&\leq \mathbb{P}_W(\phi((\overline{\mathbf{Y} - \mu})^{\langle W \rangle}) > q_{\alpha_0(1-\delta)}(\phi, \mathbf{Y} - \overline{\mathbf{Y}}) + \gamma_1 f(\mathbf{Y})) \\
&\leq \mathbb{P}_W(\phi((\overline{\mathbf{Y} - \overline{\mathbf{Y}}})^{\langle W \rangle}) > q_{\alpha_0(1-\delta)}(\phi, \mathbf{Y} - \overline{\mathbf{Y}})) \\
&\quad + \mathbb{P}_W(\phi(\overline{W}(\overline{\mathbf{Y}} - \mu)) > \gamma_1 f(\mathbf{Y})) \\
&\leq \alpha_0(1-\delta) + \mathbb{P}_W(\phi(\overline{W}(\overline{\mathbf{Y}} - \mu)) > \gamma_1 f(\mathbf{Y})).
\end{aligned}$$

For the first and last inequalities, we have used (29) and for the second inequality, the definition of $\mathcal{E}^c$. From this, we deduce that

$$\mathcal{E}^c \subset \{\mathbf{Y} | \mathbb{P}_W(\phi(\overline{W}(\overline{\mathbf{Y}} - \mu)) > \gamma_1 f(\mathbf{Y})) \geq \alpha_0 \delta\}.$$

Now, using the positive homogeneity of $\phi$ and the fact that both $\phi$ and $f$ are nonnegative, we have

$$\begin{aligned}
\mathbb{P}_W(\phi(\overline{W}(\overline{\mathbf{Y}} - \mu)) &> \gamma_1 f(\mathbf{Y})) \\
&= \mathbb{P}_W\left(|\overline{W}| > \frac{\gamma_1 f(\mathbf{Y})}{\phi(\operatorname{sign}(\overline{W})(\overline{\mathbf{Y}} - \mu))}\right) \\
&\leq \mathbb{P}_W\left(|\overline{W}| > \frac{\gamma_1 f(\mathbf{Y})}{\widetilde{\phi}(\overline{\mathbf{Y}} - \mu)}\right) \\
&= 2\mathbb{P}_{B_n}\left(\frac{1}{n}(2B_n - n) > \frac{\gamma_1 f(\mathbf{Y})}{\widetilde{\phi}(\overline{\mathbf{Y}} - \mu)}\right),
\end{aligned}$$

where $B_n$ denotes a Binomial$(n, \frac{1}{2})$ variable (independent of $\mathbf{Y}$). From the last two displays and the definition of $\gamma_1$, we conclude that $\mathcal{E}^c \subset \{\mathbf{Y} | \widetilde{\phi}(\overline{\mathbf{Y}} - \mu) > f(\mathbf{Y})\}$, which, substituted back into (31), leads to the desired conclusion. □

PROOF OF COROLLARY 3.3. Define the function

$$g_0(\mathbf{Y}) = q_{(1-\delta)\alpha_0}(\phi, \mathbf{Y} - \overline{\mathbf{Y}}) + \left(\sum_{i=1}^{J-1} \gamma_i q_{(1-\delta)\alpha_i}(\widetilde{\phi}, \mathbf{Y} - \overline{\mathbf{Y}}) + \gamma_J f(\mathbf{Y})\right)$$

and, for $k = 1, \ldots, J$,

$$g_k(\mathbf{Y}) = \gamma_k^{-1}\left(\sum_{i=k}^{J-1} \gamma_i q_{(1-\delta)\alpha_i}(\widetilde{\phi}, \mathbf{Y} - \overline{\mathbf{Y}}) + \gamma_J f(\mathbf{Y})\right)$$

with the convention that $g_J = f$. For $0 \leq k \leq J-1$, applying Theorem 3.2 with the function $g_{k+1}$ yields the relation

$$\mathbb{P}_W(\phi(\overline{\mathbf{Y}} - \mu) > g_k(\mathbf{Y})) \leq \alpha_k + \mathbb{P}_W(\phi(\overline{\mathbf{Y}} - \mu) > g_{k+1}(\mathbf{Y})).$$

Therefore, we get

$$\mathbb{P}_W(\phi(\overline{\mathbf{Y}} - \mu) > g_0(\mathbf{Y})) \leq \sum_{i=0}^{J-1} \alpha_i + \mathbb{P}(\widetilde{\phi}(\overline{\mathbf{Y}} - \mu) > f(\mathbf{Y}))$$

as announced. □

PROOF OF PROPOSITION 3.4. Let us first prove that an analog of Lemma 3.1 holds with $q_{\alpha_0}$ replaced by $\widetilde{q}_{\alpha_0}$. First, we have

$$\mathbb{E}_\mathbf{W} \mathbb{P}_\mathbf{Y}(\phi(\overline{\mathbf{Y}} - \mu) > \widetilde{q}_{\alpha_0}(\phi, \mathbf{Y} - \mu, \mathbf{W}))$$
$$= \mathbb{E}_{W'} \mathbb{E}_\mathbf{W} \mathbb{P}_\mathbf{Y}(\phi((\overline{\mathbf{Y} - \mu})^{\langle W'\rangle}) > \widetilde{q}_{\alpha_0}(\phi, (\mathbf{Y} - \mu) \bullet W', \mathbf{W}))$$
$$= \mathbb{E}_\mathbf{Y} \mathbb{P}_{\mathbf{W},W'}(\phi((\overline{\mathbf{Y} - \mu})^{\langle W'\rangle}) > \widetilde{q}_{\alpha_0}(\phi, \mathbf{Y} - \mu, W' \bullet \mathbf{W})),$$

where $W'$ denotes a Rademacher vector independent of all other random variables and $W' \bullet \mathbf{W} = \text{diag}(W') \cdot \mathbf{W}$ denotes the matrix obtained by multiplying the $i$th row of $\mathbf{W}$ by $W'_i$, $i = 1, \ldots, n$. Note that $(W', W' \bullet \mathbf{W}) \sim (W', \mathbf{W})$. Therefore, by definition of the quantile $\widetilde{q}_{\alpha_0}$, the latter quantity is equal to

$$\mathbb{E}_\mathbf{Y} \mathbb{P}_{\mathbf{W},W'}\left(\frac{1}{B}\sum_{j=1}^B \mathbb{1}\{\phi((\overline{\mathbf{Y} - \mu})^{\langle \mathbf{W}^j\rangle}) \geq \phi((\overline{\mathbf{Y} - \mu})^{\langle W'\rangle})\} \leq \alpha_0\right) \leq \frac{\lfloor B\alpha_0 \rfloor + 1}{B+1},$$

where the last step comes from Lemma 5.2 (see below).

The rest of the proof is similar to the one of Theorem 3.2, where $\mathbb{P}_W$ is replaced by the empirical distribution based on $\mathbf{W}$, $\widetilde{\mathbb{P}}_\mathbf{W} = \frac{1}{B}\sum_{j=1}^B \delta_{\mathbf{W}^j}$. Thus, (29) becomes, for any fixed $\mathbf{Y}, \mathbf{W}$,

$$\widetilde{\mathbb{P}}_\mathbf{W}[\phi(\overline{\mathbf{Y}}^{\langle W\rangle}) > \widetilde{q}_{\alpha_0}(\phi, \mathbf{Y}, \mathbf{W})] \leq \alpha_0 \leq \widetilde{\mathbb{P}}_\mathbf{W}[\phi(\overline{\mathbf{Y}}^{\langle W\rangle}) \geq \widetilde{q}_{\alpha_0}(\phi, \mathbf{Y}, \mathbf{W})].$$

The role of $\mathcal{E}$ is then taken by

$$\widetilde{\mathcal{E}} := \{\mathbf{Y}, \mathbf{W} | \widetilde{q}_{\alpha_0}(\phi, \mathbf{Y} - \mu, \mathbf{W}) \leq \widetilde{q}_{\alpha_0(1-\delta)}(\phi, \mathbf{Y} - \overline{\mathbf{Y}}, \mathbf{W}) + \gamma f(\mathbf{Y}, \mathbf{W})\},$$

where we write $\gamma = \gamma(\mathbf{W}, \alpha_0 \delta)$ for short. We then have, similarly to (31),

$$\mathbb{P}_{\mathbf{Y},\mathbf{W}}(\phi(\overline{\mathbf{Y}} - \mu) > \widetilde{q}_{\alpha_0(1-\delta)}(\phi, \mathbf{Y} - \overline{\mathbf{Y}}) + \gamma f(\mathbf{Y}, \mathbf{W})) \leq \frac{\lfloor B\alpha_0 \rfloor + 1}{B+1} + \mathbb{P}_{\mathbf{Y},\mathbf{W}}(\widetilde{\mathcal{E}}^c)$$

and follow the proof of Theorem 3.2 further, we obtain

$$\widetilde{\mathcal{E}}^c \subset \left\{ \mathbf{Y}, \mathbf{W} \Big| \widetilde{\mathbb{P}}_{\mathbf{W}}\left[|\overline{W}| > \frac{\gamma f(\mathbf{Y}, \mathbf{W})}{\widetilde{\phi}(\overline{\mathbf{Y}} - \mu)}\right] \geq \alpha_0 \delta \right\},$$

which gives the result. $\square$

We have used the following lemma which essentially reproduces Lemma 1 of [30], with a minor strengthening. While the proof was left to the reader in [30], because it was considered either elementary or common knowledge, we include a succinct proof below for completeness.

LEMMA 5.2 (Minor variation of Lemma 1 of [30]). *Let $Z_0, Z_1, \ldots, Z_B$ be exchangeable real-valued random variables. Then, for all $\alpha \in (0,1)$,*

$$\mathbb{P}\left(\frac{1}{B}\sum_{j=1}^{B} \mathbb{1}\{Z_j \geq Z_0\} \leq \alpha\right) \leq \frac{\lfloor B\alpha \rfloor + 1}{B+1} \leq \alpha + \frac{1}{B+1}.$$

*The first inequality becomes an equality if $Z_i \neq Z_j$ a.s. For example, it is the case if the $Z_i$'s are i.i.d. variables from a distribution without atoms.*

PROOF. Let $U$ denote a random variable uniformly distributed in $\{0, \ldots, B\}$ and independent of the $Z_i$'s. We then have

$$\mathbb{P}\left(\frac{1}{B}\sum_{j=1}^{B}\mathbb{1}\{Z_j \geq Z_0\} \leq \alpha\right)$$

$$= \mathbb{P}\left(\sum_{j=0}^{B}\mathbb{1}\{Z_j \geq Z_0\} \leq B\alpha + 1\right)$$

$$= \mathbb{P}_U \mathbb{P}_{(Z_i)}\left(\sum_{j=0}^{B}\mathbb{1}\{Z_j \geq Z_U\} \leq B\alpha + 1\right)$$

$$= \mathbb{P}_{(Z_i)} \mathbb{P}_U\left(\sum_{j=0}^{B}\mathbb{1}\{Z_j \geq Z_U\} \leq \lfloor B\alpha \rfloor + 1\right) \leq \frac{\lfloor B\alpha \rfloor + 1}{B+1}.$$

Note that the last inequality is an equality if the $Z_i$'s are a.s. distinct. $\square$

5.3. *Exchangeable resampling computations.* In this section, we compute constants $A_W$, $B_W$, $C_W$ and $D_W$ [defined by (3) to (6)] for some exchangeable resamplings. This implies all of the statements in Table 1. We first define several additional exchangeable resampling weights (normalized so that $\mathbb{E}[W_i] = 1$):

- *Bernoulli(p)*, $p \in (0,1)$: $pW_i$ i.i.d. with a Bernoulli distribution of parameter $p$. A classical choice is $p = \frac{1}{2}$.
- *Efron(q)*, $q \in \{1,\ldots,n\}$: $qn^{-1}W$ has a multinomial distribution with parameters $(q; n^{-1},\ldots,n^{-1})$. A classical choice is $q = n$.
- *Poisson(µ)*, $\mu \in (0,+\infty)$: $\mu W_i$ i.i.d. with a Poisson distribution of parameter $\mu$. A classical choice is $\mu = 1$.

Note that $\overline{\mathbf{Y}}^{\langle W - \overline{W}\rangle}$ and all of the resampling constants are invariant under translation of the weights so that Bernoulli(1/2) weights are completely equivalent to Rademacher weights in this paper.

LEMMA 5.3. 1. *Let $W$ be Bernoulli(p) weights with $p \in (0,1)$. We then have* $2(1-p)(1-\frac{1}{n}) = A_W \leq B_W \leq \sqrt{\frac{1}{p}-1}\sqrt{1-\frac{1}{n}}$, $C_W = \sqrt{\frac{1}{p}-1}$ *and* $D_W \leq \frac{1}{2p} + |\frac{1}{2p} - 1| + \sqrt{\frac{1-p}{np}}$.

2. *Let $W$ be Efron(q) weights with $q \in \{1,\ldots,n\}$. We then have* $2(1-\frac{1}{n})^q = A_W \leq B_W \leq \sqrt{\frac{n-1}{q}}$ *and* $C_W = \sqrt{\frac{n}{q}}$.

3. *Let $W$ be Poisson(µ) weights with $\mu > 0$. We then have* $A_W \leq B_W \leq \frac{1}{\sqrt{\mu}}\sqrt{1-\frac{1}{n}}$ *and* $C_W = \frac{1}{\sqrt{\mu}}$. *Moreover, if $\mu = 1$, we get* $\frac{2}{e} - \frac{1}{\sqrt{n}} \leq A_W$.

4. *Let $W$ be Random hold-out(q) weights with $q \in \{1,\ldots,n\}$. We then have* $A_W = 2(1-\frac{q}{n})$, $B_W = \sqrt{\frac{n}{q}-1}$, $C_W = \sqrt{\frac{n}{n-1}}\sqrt{\frac{n}{q}-1}$ *and* $D_W = \frac{n}{2q} + |1 - \frac{n}{2q}|$.

PROOF. We consider the following cases:

*General case.* First, we only assume that $W$ is exchangeable. Then, from the concavity of $\sqrt{\cdot}$ and the triangular inequality, we have

$$
\begin{aligned}
\mathbb{E}|W_1 - \mathbb{E}[W_1]| &- \sqrt{\mathbb{E}(\overline{W} - \mathbb{E}[W_1])^2} \\
&\leq \mathbb{E}|W_1 - \mathbb{E}[W_1]| - \mathbb{E}|\overline{W} - \mathbb{E}[W_1]| \leq A_W \leq B_W \leq \sqrt{\frac{n-1}{n}}C_W.
\end{aligned}
\tag{32}
$$

*Independent weights.* When we suppose that the $W_i$ are i.i.d., we get

$$\mathbb{E}|W_1 - \mathbb{E}[W_1]| - \frac{\sqrt{\mathrm{Var}(W_1)}}{\sqrt{n}} \leq A_W \quad \text{and} \quad C_W = \sqrt{\mathrm{Var}(W_1)}. \tag{33}$$

*Bernoulli.* First, we have $A_W = \mathbb{E}|W_1 - \overline{W}| = \mathbb{E}|(1 - \frac{1}{n})W_1 - X_{n,p}|$ with $X_{n,p} := \frac{1}{n}(W_2 + \cdots + W_n)$. Since $W_1$ and $X_{n,p}$ are independent and $X_{n,p} \in [0, (n-1)/(np)]$ a.s., we obtain

$$A_W = p\mathbb{E}\left[\left(1 - \frac{1}{n}\right)\frac{1}{p} - X_{n,p}\right] + (1-p)\mathbb{E}[X_{n,p}] = 1 - \frac{1}{n} + (1 - 2p)\mathbb{E}[X_{n,p}].$$

The formula for $A_W$ follows since $\mathbb{E}[X_{n,p}] = (n-1)/n$. Second, note that the Bernoulli($p$) weights are i.i.d. with $\mathrm{Var}(W_1) = p^{-1} - 1$, $\mathbb{E}[W_1] = 1$ and $\mathbb{E}|W_1 - 1| = p(p^{-1} - 1) + (1 - p) = 2(1 - p)$. Hence, (32) and (33) lead to the bounds for $B_W$ and $C_W$. Finally, the Bernoulli($p$) weights satisfy the assumption of (6) with $x_0 = a = (2p)^{-1}$. Then

$$D_W = \frac{1}{2p} + \mathbb{E}\left|\overline{W} - \frac{1}{2p}\right| \leq \frac{1}{2p} + \left|1 - \frac{1}{2p}\right| + \mathbb{E}|\overline{W} - 1|$$

$$\leq \frac{1}{2p} + \frac{1}{p}\left|\frac{1}{2} - p\right| + \sqrt{\frac{1-p}{np}}.$$

*Efron.* We have $\overline{W} = 1$ a.s. so that $C_W = \sqrt{\frac{n}{n-1} \times \mathrm{Var}(W_1)} = \sqrt{n/q}$. If, moreover, $q \leq n$, then $W_i < 1$ implies that $W_i = 0$ and $A_W = \mathbb{E}|W_1 - 1| = \mathbb{E}[W_1 - 1 + 2\mathbb{1}\{W_1 = 0\}] = 2\mathbb{P}(W_1 = 0) = 2(1 - \frac{1}{n})^q$. The result follows from (32).

*Poisson.* These weights are i.i.d. with $\mathrm{Var}(W_1) = \mu^{-1}$, $\mathbb{E}[W_1] = 1$. Moreover, if $\mu \leq 1$, $W_i < 1$ implies that $W_i = 0$ and $\mathbb{E}|W_1 - 1| = 2\mathbb{P}(W_1 = 0) = 2e^{-\mu}$. With (32) and (33), the result follows.

*Random hold-out.* These weights are such that $\{W_i\}_{1 \leq i \leq n}$ takes only two values, with $\overline{W} = 1$. Then $A_W$, $B_W$ and $C_W$ can be directly computed. Moreover, they satisfy the assumption of (6) with $x_0 = a = n/(2q)$. The computation of $D_W$ is straightforward. □

**Acknowledgments.** The first author's research was mostly carried out at University Paris-Sud (Laboratoire de Mathematiques, CNRS UMR 8628). The second author's research was partially carried out while holding an invited position at the Statistics Department of the University of Chicago, which is warmly acknowledged. The third author's research was mostly carried out at the French institute INRA-Jouy and at the Free University of Amsterdam. We wish to thank Pascal Massart for his particularly relevant comments and suggestions. We would also like to thank the two referees and the Associate Editor for their insights which led, in particular, to a more rational organization of the paper.

# REFERENCES

[1] Arlot, S. (2007). *Resampling and Model Selection*. Ph.D. thesis, Univ. Paris XI.
[2] Arlot, S., Blanchard, G. and Roquain, É. (2010). Some nonasymptotic results on resampling in high dimension. II: Multiple tests. *Ann. Statist.* **38** 83–99.
[3] Baraud, Y. (2004). Confidence balls in Gaussian regression. *Ann. Statist.* **32** 528–551. MR2060168
[4] Beran, R. (2003). The impact of the bootstrap on statistical algorithms and theory. *Statist. Sci.* **18** 175–184. MR2026078
[5] Beran, R. and Dümbgen, L. (1998). Modulation of estimators and confidence sets. *Ann. Statist.* **26** 1826–1856. MR1673280
[6] Cai, T. and Low, M. (2006). Adaptive confidence balls. *Ann. Statist.* **34** 202–228. MR2275240
[7] Cirel'son, B. R., Ibragimov, I. A. and Sudakov, V. N. (1976). Norms of Gaussian sample functions. In *Proceedings of the Third Japan–USSR Symposium on Probability Theory. Lecture Notes in Mathematics* **550** 20–41. Springer, Berlin. MR0458556
[8] Darvas, F., Rautiainen, M., Pantazis, D., Baillet, S., Benali, H., Mosher, J., Garnero, L. and Leahy, R. (2005). Investigations of dipole localization accuracy in MEG using the bootstrap. *NeuroImage* **25** 355–368.
[9] DiCiccio, T. J. and Efron, B. (1996). Bootstrap confidence intervals. *Statist. Sci.* **11** 189–228. MR1436647
[10] Durot, C. and Rozenholc, Y. (2006). An adaptive test for zero mean. *Math. Methods Statist.* **15** 26–60. MR2225429
[11] Efron, B. (1979). Bootstrap methods: Another look at the jackknife. *Ann. Statist.* **7** 1–26. MR0515681
[12] Fisher, R. A. (1935). *The Design of Experiments.* Oliver and Boyd, Edinburgh.
[13] Fromont, M. (2007). Model selection by bootstrap penalization for classification. *Mach. Learn.* **66** 165–207.
[14] Ge, Y., Dudoit, S. and Speed, T. P. (2003). Resampling-based multiple testing for microarray data analysis. *Test* **12** 1–77. MR1993286
[15] Hall, P. (1992). *The Bootstrap and Edgeworth Expansion*. Springer, New York. MR1145237
[16] Hall, P. and Mammen, E. (1994). On general resampling algorithms and their performance in distribution estimation. *Ann. Statist.* **22** 2011–2030. MR1329180
[17] Hoffmann, M. and Lepski, O. (2002). Random rates in anisotropic regression. *Ann. Statist.* **30** 325–396. MR1902892
[18] Jerbi, K., Lachaux, J.-P., N'Diaye, K., Pantazis, D., Leahy, R. M., Garnero, L. and Baillet, S. (2007). Coherent neural representation of hand speed in humans revealed by MEG imaging. *PNAS* **104** 7676–7681.
[19] Juditsky, A. and Lambert-Lacroix, S. (2003). Nonparametric confidence set estimation. *Math. Methods Statist.* **12** 410–428. MR2054156
[20] Koltchinskii, V. (2001). Rademacher penalties and structural risk minimization. *IEEE Trans. Inform. Theory* **47** 1902–1914. MR1842526
[21] Lepski, O. V. (1999). How to improve the accuracy of estimation. *Math. Methods Statist.* **8** 441–486. MR1755896
[22] Li, K.-C. (1989). Honest confidence regions for nonparametric regression. *Ann. Statist.* **17** 1001–1008. MR1015135
[23] Mason, D. M. and Newton, M. A. (1992). A rank statistics approach to the consistency of a general bootstrap. *Ann. Statist.* **20** 1611–1624. MR1186268

[24] MASSART, P. (2007). *Concentration Inequalities and Model Selection (Lecture Notes of the St-Flour Probability Summer School 2003). Lecture Notes in Mathematics* **1896**. Springer, Berlin. MR2319879
[25] MCDIARMID, C. (1989). On the method of bounded differences. In *Surveys in Combinatorics. London Mathematical Society Lecture Notes* **141** 148–188. Cambridge Univ. Press, Cambridge. MR1036755
[26] PANTAZIS, D., NICHOLS, T. E., BAILLET, S. and LEAHY, R. M. (2005). A comparison of random field theory and permutation methods for statistical analysis of MEG data. *NeuroImage* **25** 383–394.
[27] POLITIS, D. N., ROMANO, J. P. and WOLF, M. (1999). *Subsampling*. Springer, New York. MR1707286
[28] PRÆSTGAARD, J. and WELLNER, J. A. (1993). Exchangeably weighted bootstraps of the general empirical process. *Ann. Probab.* **21** 2053–2086. MR1245301
[29] ROBINS, J. and VAN DER VAART, A. (2006). Adaptive nonparametric confidence sets. *Ann. Statist.* **34** 229–253. MR2275241
[30] ROMANO, J. P. and WOLF, M. (2005). Exact and approximate stepdown methods for multiple hypothesis testing. *J. Amer. Statist. Assoc.* **100** 94–108. MR2156821
[31] VAN DER VAART, A. W. and WELLNER, J. A. (1996). *Weak Convergence and Empirical Processes*. Springer, New York. MR1385671
[32] WABERSKI, T., GOBBELE, R., KAWOHL, W., CORDES, C. and BUCHNER, H. (2003). Immediate cortical reorganization after local anesthetic block of the thumb: Source localization of somatosensory evoked potentials in human subjects. *Neurosci. Lett.* **347** 151–154.

S. ARLOT
CNRS: WILLOW PROJECT-TEAM
LABORATOIRE D'INFORMATIQUE
  DE L'ECOLE NORMALE SUPERIEURE
(CNRS/ENS/INRIA UMR 8548)
INRIA, 23 AVENUE D'ITALIE, CS 81321
75214 PARIS CEDEX 13
FRANCE
E-MAIL: sylvain.arlot@ens.fr

G. BLANCHARD
WEIERSTRASS INSTITUTE
  FOR APPLIED STOCHASTICS
  AND ANALYSIS
MOHRENSTRASSE 39, 10117 BERLIN
GERMANY
E-MAIL: blanchar@wias-berlin.de

E. ROQUAIN
UPMC UNIVERSITY OF PARIS 6
UMR 7599, LPMA
4, PLACE JUSSIEU
75252 PARIS CEDEX 05
FRANCE
E-MAIL: etienne.roquain@upmc.fr

# SOME NONASYMPTOTIC RESULTS ON RESAMPLING IN HIGH DIMENSION, II: MULTIPLE TESTS[1]

By Sylvain Arlot, Gilles Blanchard[2] and Etienne Roquain

*CNRS ENS, Weierstrass Institute Berlin and UPMC University of Paris 6*

In the context of correlated multiple tests, we aim to nonasymptotically control the family-wise error rate (FWER) using resampling-type procedures. We observe repeated realizations of a Gaussian random vector in possibly high dimension and with an unknown covariance matrix, and consider the one- and two-sided multiple testing problem for the mean values of its coordinates. We address this problem by using the confidence regions developed in the companion paper [*Ann. Statist.* (2009), to appear], which lead directly to single-step procedures; these can then be improved using step-down algorithms, following an established general methodology laid down by Romano and Wolf [*J. Amer. Statist. Assoc.* **100** (2005) 94–108]. This gives rise to several different procedures, whose performances are compared using simulated data.

## 1. Introduction.

1.1. *Framework and motivations.* We consider a sample $\mathbf{Y} := (\mathbf{Y}^1, \ldots, \mathbf{Y}^n)$ of $n \geq 2$ i.i.d. observations of a Gaussian vector with dimensionality $K$, possibly much larger than $n$. The common covariance matrix of the $\mathbf{Y}^i$ is not assumed to be known in advance. We investigate the two following multiple testing problems for the common mean $\mu \in \mathbb{R}^K$ of the $\mathbf{Y}^i$:

- *One-sided.* Test simultaneously $H_k$: "$\mu_k \leq 0$" against $A_k$: "$\mu_k > 0$" for $1 \leq k \leq K$;

Received November 2007; revised August 2008.
[1]Supported in part by the IST and ICT programs of the European Community, respectively, under the PASCAL (IST-2002-506778) and PASCAL2 (ICT-216886) Networks of Excellence.
[2]Supported in part by the Fraunhofer Institute First, Berlin.
*AMS 2000 subject classifications.* Primary 62G10; secondary 62G09.
*Key words and phrases.* Family-wise error, multiple testing, high-dimensional data, nonasymptotic error control, resampling, resampled quantile.

- *Two-sided.* Test simultaneously $H_k$: "$\mu_k = 0$" against $A_k$: "$\mu_k \neq 0$" for $1 \leq k \leq K$.

For simplicity, we introduce the following notation to cover both cases:

(1) $\quad$ test simultaneously $H_k$: "$[\![\mu_k]\!] = 0$" against $A_k$: "$[\![\mu_k]\!] \neq 0$" for $1 \leq k \leq K$,

where, for $x \in \mathbb{R}$, $[\![x]\!]$ denotes either $\max\{x, 0\} = x_+$ in the one-sided context or $|x|$ in the two-sided context.

In this paper, we tackle the problem (1) by building multiple testing procedures which control the family-wise error rate (FWER). We emphasize that:

- we aim to obtain a *nonasymptotic* control, valid for any fixed $K$ and $n$, and, in particular, with $K$ possibly much larger than the number of observations $n$;
- we do not want to make any particular prior assumption on the structure of the covariance matrix of the $\mathbf{Y}^i$.

As explained in [1], this point of view is motivated by some practical applications, especially neuroimaging [5, 10, 11]. Multiple testing problems in this field typically have parameters $10^4 \leq K \leq 10^7$, $n \leq 100$, with strong and complex dependencies between the coordinates of $\mathbf{Y}^i$. Another motivating example is microarray data analysis (see, e.g., [8]).

1.2. *Goals.* In this work, we consider thresholding-based procedures which reject the null hypotheses $H_k$ for indices $k \in R_\alpha(\mathbf{Y}) \subset \mathcal{H} := \{1, \ldots, K\}$ corresponding to large values of $[\![\overline{\mathbf{Y}}_k]\!]$, where $\overline{\mathbf{Y}}_k = n^{-1} \sum_{i=1}^n \mathbf{Y}_k^i$ denotes the vector of empirical means, that is,

(2) $\quad R_\alpha(\mathbf{Y}) = \{1 \leq k \leq K | [\![\overline{\mathbf{Y}}_k]\!] > t_\alpha(\mathbf{Y})\}$,

where $t_\alpha(\mathbf{Y})$ is a possibly data-dependent threshold.

The type I error of such a multiple testing procedure is measured here by the family-wise error rate (FWER), defined as the probability that at least one hypothesis is wrongly rejected:

$$\mathrm{FWER}(R_\alpha) := \mathbb{P}(R_\alpha(\mathbf{Y}) \cap \mathcal{H}_0 \neq \varnothing),$$

where $\mathcal{H}_0 := \{k | [\![\mu_k]\!] = 0\}$ is the set of coordinates corresponding to the true null hypotheses. The choice of this error rate is discussed in Section 5.1. Given a level $\alpha \in (0, 1)$, the goal is now to build a multiple testing procedure $R_\alpha$ such that $\mathrm{FWER}(R_\alpha) \leq \alpha$ is valid for all distributions in the family being considered (i.e., Gaussian with arbitrary mean vector and covariance matrix); furthermore, as many false hypotheses as possible should be rejected.

To this end, we use the family of $(1 - \alpha)$-resampling-based confidence regions for $\mu$ introduced in the companion paper [1]. Of interest here are regions taking the following form: for some subset $\mathcal{C} \subset \mathcal{H}$,

$$(3) \qquad \mathcal{G}(\mathbf{Y}, 1 - \alpha, \mathcal{C}) := \Big\{ x \in \mathbb{R}^K | \sup_{k \in \mathcal{C}} [\![\overline{\mathbf{Y}}_k - x_k]\!] \leq t_\alpha(\mathbf{Y}, \mathcal{C}) \Big\},$$

where $t_\alpha$ is a data-dependent threshold built using a resampling principle. Several possible choices for this threshold were proposed. The main results of [1], as well as the link between confidence regions (3) and (single-step) multiple tests for (1), are briefly recalled in Section 2.

1.3. *Contribution in relation to previous work.* Most of the existing resampling-type multiple testing procedures have been developed in an asymptotic framework (see, e.g., [8, 13, 17–19]), while our present goal is to study procedures that have nonasymptotic theoretical validity (for any $K$ and $n$). The main classical alternative approach to asymptotic validity is to use an invariance of the null distribution under a group of transformations—that is, exact randomized tests [14–16] (the underlying idea can be traced back to Fisher's permutation test [7]). Additionally, and as explained in [16], exact tests can be combined with a step-down algorithm to build less conservative procedures while preserving the same nonasymptotic control on their FWER (also, see [17] for a generalization to the $k$-FWER).

In the case considered here, Gaussian vectors $\mathbf{Y}^i$ have a symmetric distribution around their mean so that the action of mirroring any subset of the vectors in the data sample with respect to their mean constitutes such a group of distribution-preserving transformations. In the two-sided case, this group is known under the global null hypothesis $\mu = 0$ and just corresponds to arbitrary sign reversal of each data vector.

Consequently, it is possible to directly derive from [16] a step-down procedure whose FWER is controlled in a nonasymptotic setting (see Section 3). This approach will be referred to as *uncentered* in this paper because the sign reversal is applied to the $(\mathbf{Y}^i)_{1 \leq i \leq n}$ themselves, without prior centering.

We observe that the principle of sign reversal was also used in [6] in order to build an adaptive (single) test for zero mean under the assumption of symmetric and independent errors. The setting studied here is different since we consider multiple testing with possibly dependent errors.

Compared to this uncentered approach, most of the procedures proposed in this paper consist of applying the sign reversal to the *empirically centered data* $(\mathbf{Y}^i - \overline{\mathbf{Y}})_{1 \leq i \leq n}$. It was proven in [1] that such an intuitive idea is theoretically valid, despite the dependencies between the $\mathbf{Y}^i - \overline{\mathbf{Y}}$, $1 \leq i \leq n$, at the cost of adding a second order remainder term. We argue in the present paper that in some interesting situations, the prior centering operation leads to a noticeable decrease of the computation time of the step-down algorithm, up to some small loss in accuracy (due to the remainder term) with respect to the uncentered step-down. Additionally, the centered approach can be used both in the one-sided and two-sided contexts, while the uncentered approach has, to the best of our knowledge, only been proven valid in the two-sided case.

1.4. *Notation.* Let us now introduce some notation that will be used throughout this paper.

- $\mathbf{Y}$ denotes the $K \times n$ data matrix $(\mathbf{Y}_k^i)_{1 \leq k \leq K, 1 \leq i \leq n}$. A superscript index such as $\mathbf{Y}^i$ indicates the $i$th column of a matrix. The empirical mean vector is $\overline{\mathbf{Y}} := \frac{1}{n} \sum_{i=1}^n \mathbf{Y}^i$. If $\mu \in \mathbb{R}^K$, $\mathbf{Y} - \mu$ is the matrix obtained by subtracting $\mu$ from each (column) vector of $\mathbf{Y}$.
- The vector $\sigma := (\sigma_k)_{1 \leq k \leq K}$ is the vector of the standard deviations of the data: $\forall k, 1 \leq k \leq K$, $\sigma_k := \mathrm{Var}^{1/2}(\mathbf{Y}_k^1)$. For $\mathcal{C} \subset \mathcal{H}$, we also define $\|\sigma\|_\mathcal{C} := \sup_{k \in \mathcal{C}} \sigma_k$.
- $\overline{\Phi}$ is the standard Gaussian upper tail function: if $X \sim \mathcal{N}(0, 1)$, then $\forall x \in \mathbb{R}$, $\overline{\Phi}(x) = \mathbb{P}(X \geq x)$.
- If $W \in \mathbb{R}^n$, we define the mean of $W \in \mathbb{R}^n$ as $\overline{W} := \frac{1}{n} \sum_{i=1}^n W_i$ and for every $c \in \mathbb{R}$, $W - c := (W_i - c)_{1 \leq i \leq n} \in \mathbb{R}^n$.
- For a subset $\mathcal{C} \subset \mathcal{H}$, $|\mathcal{C}|$ denotes the cardinality of $\mathcal{C}$.

## 2. Single-step procedures using resampling-based thresholds.

2.1. *Connection between confidence regions and FWER control.* We start with recalling a simple device linking confidence regions to FWER control in multiple testing. In a nutshell, the idea is that a confidence region of the form (3) directly gives a multiple testing procedure $R$ with controlled FWER when taking $\mathcal{C} = \mathcal{H}_0$. Since $\mathcal{H}_0$ is not known in advance, we actually need a confidence region (3) defined for every $\mathcal{C} \subset \mathcal{H}$ and satisfying certain properties.

More formally, let $\alpha \in (0, 1)$ be fixed and $\mathcal{T}_\alpha = (t_\alpha(\mathbf{Y}, \mathcal{C}), \mathcal{C} \subset \mathcal{H}, \mathbf{Y} \in \mathbb{R}^{K \times n})$ be a family of thresholds indexed by subsets $\mathcal{C} \subset \mathcal{H}$. We consider threshold families satisfying the two following key properties. First, $t_\alpha(\mathbf{Y}, \mathcal{H}_0)$ is a $1 - \alpha$ confidence bound on the deviations of $\sup_{k \in \mathcal{H}_0} [\overline{\mathbf{Y}}_k]$:

$$(\mathrm{CB}_\alpha) \qquad \mathbb{P}\Big(\sup_{k \in \mathcal{H}_0} [\overline{\mathbf{Y}}_k] < t_\alpha(\mathbf{Y}, \mathcal{H}_0)\Big) \geq 1 - \alpha.$$

Second, $\mathcal{T}_\alpha$ is nondecreasing w.r.t. $\mathcal{C}$, that is,

$$(\mathrm{ND}) \quad \forall \mathbf{Y} \in \mathbb{R}^K, \forall \mathcal{C} \quad \mathcal{C}' \subset \mathcal{H}, \mathcal{C} \subset \mathcal{C}' \quad \Rightarrow \quad t_\alpha(\mathbf{Y}, \mathcal{C}) \leq t_\alpha(\mathbf{Y}, \mathcal{C}').$$

We now define a single-step multiple testing procedure and establish its FWER control.

PROPOSITION 2.1. *Define the single-step multiple testing procedure associated with $\mathcal{T}_\alpha$ as the procedure rejecting the set of hypotheses given by*

(4) $$\{k \in \mathcal{H} | [\![\overline{\mathbf{Y}}_k]\!] > t_\alpha(\mathbf{Y}, \mathcal{H})\}.$$

*If the threshold family satisfies ($CB_\alpha$) and (ND), then the FWER of the associated single-step procedure is controlled at level $\alpha$.*

PROOF. We first use (ND), then ($CB_\alpha$):

$$\mathbb{P}(\exists k | [\![\overline{\mathbf{Y}}_k]\!] > t_\alpha(\mathbf{Y}, \mathcal{H}) \text{ and } [\![\mu_k]\!] = 0)$$
$$= \mathbb{P}\Big(\sup_{k \in \mathcal{H}_0} [\![\overline{\mathbf{Y}}_k]\!] > t_\alpha(\mathbf{Y}, \mathcal{H})\Big)$$
$$\leq \mathbb{P}\Big(\sup_{k \in \mathcal{H}_0} [\![\overline{\mathbf{Y}}_k]\!] > t_\alpha(\mathbf{Y}, \mathcal{H}_0)\Big) \leq \alpha. \qquad \square$$

Note that the single-step procedure only uses the value of the largest threshold among the $t_\alpha(\mathbf{Y}, \mathcal{C}), \mathcal{C} \subset \mathcal{H}$. In Section 3, we use the iterative step-down principle to improve the procedure by making use of the thresholds $t_\alpha(\mathbf{Y}, \mathcal{C})$ for some smaller $\mathcal{C} \subset \mathcal{H}$.

The condition ($CB_\alpha$) is, in particular, satisfied whenever, for any $\mathcal{C} \subset \mathcal{H}$, $t(\mathbf{Y}, \mathcal{C})$ provides a $1 - \alpha$ confidence region of the form (3) for $(\mu_k)_{k \in \mathcal{C}}$. We use this idea next to derive testing thresholds from the confidence regions constructed in [1].

2.2. *Resampling thresholds.* We first give a compact recapitulation of resampling-based thresholds introduced in [1] and used to build confidence regions for the mean of a high-dimensional, correlated vector. This is intended as a single overall reference for all of the thresholds that we use in the present paper. Here, and in the following, $W \in \mathbb{R}^n$ denotes a random vector independent from the data $\mathbf{Y}$, called the *resampling weight vector*. Moreover, in order to simplify the results of [1], we specifically assume that the $W_i$ are i.i.d. Rademacher random variables, that is, that they satisfy $\mathbb{P}(W_i = 1) = \mathbb{P}(W_i = -1) = 1/2$. As first building blocks, define the two following resampling quantities, the (scaled) resampled expectation and quantile:

(5) $$\mathcal{E}(\mathbf{Y}, \mathcal{C}) := B_W^{-1} \mathbb{E}_W \left[ \sup_{k \in \mathcal{C}} \left[\!\!\left[ \frac{1}{n} \sum_{i=1}^n W_i \mathbf{Y}_k^i \right]\!\!\right] \right],$$

(6) $$q_\alpha(\mathbf{Y}, \mathcal{C}) := \inf \left\{ x \in \mathbb{R} \bigg| \mathbb{P}_W \left( \sup_{k \in \mathcal{C}} \left[\!\!\left[ \frac{1}{n} \sum_{i=1}^n W_i \mathbf{Y}_k^i \right]\!\!\right] > x \right) \leq \alpha \right\},$$

where $B_W := \mathbb{E}_W[(\frac{1}{n} \sum_{i=1}^n (W_i - \overline{W}))^{1/2}]$ and $\mathbb{E}_W[\cdot]$ [resp., $\mathbb{P}_W(\cdot)$] denotes the expectation (resp., probability) operator over the distribution of $W$ only. We

TABLE 1
*Reference table for the different rejection thresholds*

$$
(7) \quad t_{\alpha,\mathrm{Bonf}}(\mathbf{Y},\mathcal{C}) := \frac{1}{\sqrt{n}}\|\sigma\|_{\mathcal{C}}\overline{\Phi}^{-1}\left(\frac{\alpha}{c|\mathcal{C}|}\right) \quad \text{with} \quad \begin{cases} c=1 & \text{(one-sided case)} \\ c=2 & \text{(two-sided case)} \end{cases}
$$

$$
(8) \quad t_{\alpha,\mathrm{conc}}(\mathbf{Y},\mathcal{C}) := \mathcal{E}(\mathbf{Y}-\overline{\mathbf{Y}},\mathcal{C}) + \|\sigma\|_{\mathcal{C}}\overline{\Phi}^{-1}(\alpha/2)\left[\frac{1}{nB_W} + \frac{1}{\sqrt{n}}\right]
$$

$$
(9) \quad t_{\alpha,\mathrm{conc}\wedge\mathrm{Bonf}}(\mathbf{Y},\mathcal{C}) := \min\Big(t_{\alpha(1-\delta),\mathrm{Bonf}}(\mathbf{Y},\mathcal{C}),
$$
$$
\mathcal{E}(\mathbf{Y}-\overline{\mathbf{Y}},\mathcal{C}) + \frac{\|\sigma\|_{\mathcal{C}}}{\sqrt{n}}\overline{\Phi}^{-1}\left(\frac{\alpha(1-\delta)}{2}\right) + \frac{\|\sigma\|_{\mathcal{C}}}{nB_W}\overline{\Phi}^{-1}\left(\frac{\alpha\delta}{2}\right)\Big)
$$

$$
(10) \quad t^*_{\alpha,\mathrm{quant}}(\mathbf{Y},\mathcal{C}) := q_\alpha(\mathbf{Y}-\overline{\mathbf{Y}},\mathcal{C})
$$

$$
(11) \quad t_{\alpha,\mathrm{quant}+\mathrm{Bonf}}(\mathbf{Y},\mathcal{C}) := t^*_{\alpha_0(1-\delta),\mathrm{quant}}(\mathbf{Y},\mathcal{C}) + \gamma_n(\alpha_0\delta)t_{\alpha-\alpha_0,\mathrm{Bonf}}(\mathbf{Y},\mathcal{C})
$$

$$
(12) \quad t_{\alpha,\mathrm{quant}+\mathrm{conc}}(\mathbf{Y},\mathcal{C}) := t^*_{\alpha_0(1-\delta),\mathrm{quant}}(\mathbf{Y},\mathcal{C}) + \gamma_n(\alpha_0\delta)t_{\alpha-\alpha_0,\mathrm{conc}}(\mathbf{Y},\mathcal{C})
$$

$$
(13) \quad t_{\alpha,\mathrm{quant.uncent}}(\mathbf{Y},\mathcal{C}) := q_\alpha(\mathbf{Y},\mathcal{C})
$$

also define the following function which is the upper quantile function of a binomial $(n, \frac{1}{2})$ variable:

$$
\overline{\mathcal{B}}(n,\eta) := \max\left\{k \in \{0,\ldots,n\} \Big| 2^{-n}\sum_{i=k}^{n}\binom{n}{i} \geq \eta\right\}.
$$

Finally, we define the factor

$$
\gamma_n(\eta) := \frac{2\overline{\mathcal{B}}(n,\eta/2) - n}{n} \leq \left(\frac{2\log(2/\eta)}{n}\right)^{1/2},
$$

where the last inequality, intended as a more explicit formula, is obtained via Hoeffding's inequality.

Table 1 gives a reference for the different rejection thresholds considered in this paper, depending on a target type I error level $\alpha$, subset of coordinates $\mathcal{C}$ and possibly on two arbitrary parameters $\alpha_0 \in (0,\alpha)$ and $\delta \in (0,1)$. The threshold (7) is Bonferroni's for Gaussian variables. Thresholds (8), (9), (11) and (12) were introduced in [1]. More precisely, threshold (8) is based on a Gaussian concentration result. Threshold (9) is a compound threshold which is very close to the minimum of (7) and (8). Threshold (10) is a raw resampled quantile for the empirically centered data; it has not been proven theoretically that this threshold achieves the correct level (this is signalled by the star symbol). The thresholds (11) and (12) are based on the latter with an additional term which was introduced in [1] to compensate (from a theoretical point of view) for the optimism in centering the data empirically rather than using the (unknown) true mean. The thresholds (7), (8) and (9) [and thus (11) and (12)] depend on the quantity $\|\sigma\|_\mathcal{C}$; if it is unknown, a confidence upper bound on $\|\sigma\|_\mathcal{C}$ can be built (see Section 4.1 of [1]). Finally, note that all of these thresholds are nondecreasing w.r.t. $\mathcal{C}$, that is, they satisfy assumption (ND). The nonasymptotic theoretical results obtained in [1] in the Gaussian case can be summed up in the following theorem.

THEOREM 2.2. *If $t_\alpha(\mathbf{Y},\mathcal{C})$ is one of the thresholds defined by (7), (8), (9), (11) or (12), then it holds for any $\mathcal{C} \subset \mathcal{H}$, in the one-sided as well as the two-sided setting, that*

$$(14) \qquad \mathbb{P}\Big(\sup_{k\in\mathcal{C}} [\![\overline{\mathbf{Y}}_k - \mu_k]\!] < t_\alpha(\mathbf{Y},\mathcal{C})\Big) \geq 1 - \alpha.$$

*In particular, all of these thresholds satisfy ($\mathrm{CB}_\alpha$), both in the one-sided and two-sided cases.*

Note that the results obtained in [1] have more generality. In particular, variations on the above thresholds were proposed for non-Gaussian, but bounded, data and weight families different from Rademacher weights can be used in (8) and (9). For the purposes of the present work, we restrict our attention to Gaussian data and Rademacher weights for simplicity. It is straightforward to show that (14) implies ($\mathrm{CB}_\alpha$): the two-sided case is obvious since $\mu_k = 0$ for $k \in \mathcal{H}_0$; the one-sided case is an easy consequence of the fact that the positive part is a nondecreasing function. Therefore, by an application of Proposition 2.1, the corresponding thresholds $t_\alpha(\mathbf{Y},\mathcal{H})$ for the full set of hypotheses can be used for multiple testing in the one-sided as well as two-sided setting with a nonasymptotic control of the FWER.

We mentioned above that the thresholds (11) and (12), based on a resampled quantile for the *empirically centered* data $(\mathbf{Y} - \overline{\mathbf{Y}})$, include an additional term in order to upper bound the variations introduced by the centering operation. In the context of testing, however, it is important to note that the quantile for the *uncentered* data defined in (13) is (without modification) a valid threshold in the two-sided setting.

THEOREM 2.3. *Assume only that $\mathbf{Y}$ has a symmetric distribution around its mean $\mu$, that is, that $(\mathbf{Y}^1 - \mu) \sim (\mu - \mathbf{Y}^1)$. If $\mu_k = 0$ for all $k \in \mathcal{C}$, then the threshold $t_{\alpha,\mathrm{quant.uncent}}(\mathbf{Y},\mathcal{C})$ defined by (13) satisfies (14). In particular, the threshold $t_{\alpha,\mathrm{quant.uncent}}(\mathbf{Y},\mathcal{C})$ satisfies ($\mathrm{CB}_\alpha$) in the two-sided setting.*

This result can probably be considered to be well known and corresponds, for example, to Lemma 3.1 in [1]. Again by Proposition 2.1, the threshold defined by (13) can therefore be used for multiple testing (although only for the two-sided setting).

It is useful at this point to carry out a brief qualitative comparison of the uncentered quantile threshold (13) and the centered quantile thresholds (11) and (12) (in the two-sided setting). The obvious differences between the two types of thresholds are:

- the data vectors are not centered around the empirical mean $\overline{\mathbf{Y}}$ prior to computing the threshold (13);
- the centered thresholds (11) and (12) have an additional additive term with respect to the main resampled quantile; furthermore, the main centered quantile is computed at a shrunk error level $\alpha_0(1-\delta) < \alpha$.

The second point is a net drawback of the "centered" family compared to the "uncentered" one. On the other hand, empirical centering of the data has the advantage of making the corresponding threshold $t_\alpha(\mathbf{Y}, \mathcal{C})$ *translation invariant*, that is, for every $\mathbf{Y} \in \mathbb{R}^{K \times n}$ and $x \in \mathbb{R}^K$, the following property holds:

(TI) $$\forall \mathcal{C} \subset \mathcal{H} \qquad t_\alpha(\mathbf{Y} + x, \mathcal{C}) = t_\alpha(\mathbf{Y}, \mathcal{C}).$$

This property is also shared by the concentration-based thresholds (8) and (9). Therefore, large values of nonzero means $\mu_k$ do not affect these thresholds. To understand the practical consequences of this point, let us consider the following informal and qualitative argumentation. If some coordinates of $(\mathbf{Y}_k^1)_{k \in \mathcal{C}}$ have a large mean relative to the noise (i.e., a large signal-to-noise ratio or SNR), then the corresponding coordinates of $\overline{\mathbf{Y}}$ will have, on average, a large absolute value relative to the coordinates with zero mean and the contribution of the former to the threshold will make the uncentered quantile significantly larger. In contrast, the centered quantile threshold is translation invariant and thus unaffected by the signal itself. Hence, in this situation, it is likely that the centered quantile threshold will be smaller. This effect is illustrated in the simulations presented in Section 4.

**3. Step-down procedures.** Single-step procedures can often be improved by iteration based on the *step-down* principle. Roughly, the idea is to repeat the multiple testing procedure with $\mathcal{H}$ replaced by $\mathcal{H} \setminus R_\alpha(\mathbf{Y})$ and to iterate this process as long as new coordinates are rejected. Again, consider a threshold family $\mathcal{T}_\alpha = (t_\alpha(\mathbf{Y}, \mathcal{C}), \mathcal{C} \subset \mathcal{H}, \mathbf{Y} \in \mathbb{R}^{K \times n})$ satisfying (CB$_\alpha$) and (ND).

DEFINITION 3.1. Consider the nonincreasing sequence $(\mathcal{C}_j, j \geq 0)$ of subsets of $\mathcal{H}$ defined by

$$\mathcal{C}_0 := \mathcal{H} \quad \text{and} \quad \forall j \geq 1 \qquad \mathcal{C}_j := \{k \in \mathcal{C}_{j-1} | [\![\overline{\mathbf{Y}}_k]\!] \leq t_\alpha(\mathbf{Y}, \mathcal{C}_{j-1})\},$$

and let $\hat{\ell}$ be the stopping rule $\hat{\ell} = \min\{j \geq 1 | \mathcal{C}_j = \mathcal{C}_{j-1}\}$. Then the *step-down multiple testing procedure* associated with $\mathcal{T}_\alpha$ rejects the hypotheses of the set $\mathcal{H} \setminus \mathcal{C}_{\hat{\ell}}$, that is,

$$\{k \in \mathcal{H} | [\![\overline{\mathbf{Y}}_k]\!] > t_\alpha(\mathbf{Y}, \mathcal{C}_{\hat{\ell}})\}. \tag{15}$$

A very general result on step-down procedures was established in [16], Theorem 3, which we reproduce here using our notation.

THEOREM 3.2 (Romano and Wolf [16]). *Let $\mathcal{T}_\alpha$ be a threshold family satisfying (*ND*). The FWER of the step-down procedure (15) is then controlled by*

$$\mathbb{P}\Big(\sup_{k \in \mathcal{H}_0} [\![\overline{\mathbf{Y}}_k]\!] > t_\alpha(\mathbf{Y}, \mathcal{H}_0)\Big).$$

*Therefore, if $\mathcal{T}_\alpha$ additionally satisfies (*$\mathrm{CB}_\alpha$*), the FWER of the associated step-down procedure is upper bounded by $\alpha$.*

A sketch of the proof can be given as follows: assume that $\mathbf{Y}$ is such that $\sup_{k \in \mathcal{H}_0} [\![\overline{\mathbf{Y}}_k]\!] \leq t_\alpha(\mathbf{Y}, \mathcal{H}_0)$. Then $\mathcal{H}_0 \subset \mathcal{C}_{j-1}$ implies that $t_\alpha(\mathbf{Y}, \mathcal{C}_{j-1}) \geq t_\alpha(\mathbf{Y}, \mathcal{H}_0) \geq \sup_{k \in \mathcal{H}_0} [\![\overline{\mathbf{Y}}_k]\!]$ and, in turn, $\mathcal{H}_0 \subset \mathcal{C}_j$, by definition of $\mathcal{C}_j$. By recursion, $\mathcal{H}_0$ is contained in $\mathcal{C}_j$ for every $j$ and the step-down procedure therefore makes no type I error on the event $\{\sup_{k \in \mathcal{H}_0} [\![\overline{\mathbf{Y}}_k]\!] \leq t_\alpha(\mathbf{Y}, \mathcal{H}_0)\}$.

A direct consequence of Theorem 3.2 is that the step-down procedures based on any of the thresholds considered in the previous section [defined by (7), (8), (9), (11), (12) or (13)] have a FWER controlled at level $\alpha$ [for (13), only in the two-sided setting]. Note that the step-down procedure based on Bonferroni's threshold (7) is exactly Holm's procedure [9].

Parallel to the discussion at the end of Section 2.2, we can carry out a short qualitative comparison of the step-down procedure based on the uncentered quantile threshold (13) and the step-down procedures based on the centered quantile thresholds (11) and (12) (in the two-sided setting). Again, if some coordinates have a large SNR, then they certainly contribute to making the uncentered quantile threshold significantly larger at the first step of the step-down procedure. This time, however, even if this first threshold is relatively large, it will still be able to rule out at the first step precisely those coordinates having the largest means. This, in turn, will result in an important improvement of the uncentered threshold at the second iteration, and so on, until all coordinates with a large SNR have been weeded out. Thus, the initial disadvantage of the uncentered threshold will be automatically corrected along the step-down iterations and the final threshold will be close to the ideal resampled quantile $q_\alpha(\mathbf{Y}, \mathcal{H}_0)$ in the last iterations. In contrast, the centered thresholds (11) and (12) still suffer from the loss due to the remainder term and level shrinkage along the step-down. In conclusion, in contrast to the single-step case, we expect the uncentered procedure to eventually outperform the centered ones after some step-down iterations. This is in accordance with the simulations of Section 4.

At this point, it could seem that the uncentered step-down procedure is both simpler and more effective than the centered step-down ones and thus should always be preferred. However, the above discussion gives us another insight: the step-down procedure based on the uncentered quantile should require more iterations to converge because the first iterations return inaccurate thresholds. In order to deal with this drawback, we propose using the leverage of the centered quantile thresholds for the first step—weeding out in a single step most of coordinates having a large SNR—and then subsequently continuing with the uncentered threshold in the next steps for greater accuracy. We obtain the following new algorithm.

ALGORITHM 3.3 (Hybrid approach).

1. Compute the threshold $t_{\alpha,\text{quant}+\text{Bonf}}(\mathbf{Y}, \mathcal{H})$ defined by (11) with a given $\delta \in (0,1)$, $\alpha_0 \in (0, \alpha)$ and consider $R_0$, the corresponding single-step procedure (4).
2. If $R_0 = \mathcal{H}$, then stop and reject all of the null hypotheses. Otherwise, consider the set of remaining coordinates $\mathcal{C}_0 = \mathcal{H} \setminus R_0$ and apply to it the step-down procedure associated with the threshold $t_{\alpha_0,\text{quant.uncent}}(\mathbf{Y}, \mathcal{C})$ defined by (13) (at level $\alpha_0$).

PROPOSITION 3.4. *Fix $\delta \in (0,1)$ and $\alpha_0 \in (0, \alpha)$. In the two-sided context, Algorithm 3.3 gives rise to a multiple testing procedure with a FWER upper bounded by $\alpha$.*

Proposition 3.4 is proved in Section 6. What we expect is that Algorithm 3.3 essentially yields the same final result as the step-down procedure using the uncentered quantile (up to some negligible loss in the level by taking $\alpha_0$ close to $\alpha$), while requiring less iterations. In applications such as neuroimaging, where a single iteration can take up to one day, this can result in a significant improvement.

**4. Simulation study.** The (MATLAB) code used to perform the simulations of this section is available on the first author's webpage (currently at <http://www.di.ens.fr/~arlot/code/CRMTR.htm>).

4.1. *Framework.* We consider data of the form $\mathbf{Y}_t = \mu_t + G_t$, where $t$ belongs to a $d \times d$ discretized two-dimensional torus of $K = d^2$ pixels, identified with $\mathbb{T}_d^2 = (\mathbb{Z}/d\mathbb{Z})^2$, and $G$ is a centered Gaussian vector obtained by two-dimensional discrete convolution of an i.i.d. standard Gaussian field (white noise) on $\mathbb{T}_d^2$ with a function $F : \mathbb{T}_d^2 \to \mathbb{R}$ such that $\sum_{t \in \mathbb{T}_d^2} F^2(t) = 1$. This ensures that $G$ is a stationary Gaussian process on the discrete torus; it is, in particular, isotropic with $\mathbb{E}[G_t^2] = 1$ for all $t \in \mathbb{T}_d^2$.

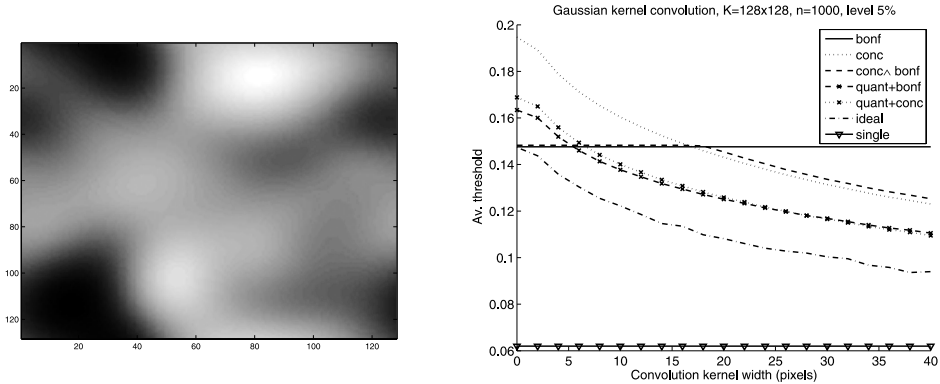

FIG. 1. *Left: example of a* $128 \times 128$ *pixel image obtained by convolution of Gaussian white noise with a pseudo-Gaussian filter with width* $b = 18$ *pixels. Right: average thresholds obtained for the different approaches; see text.*

In the simulations below, we consider, for the function $F$, a pseudo-Gaussian convolution filter of bandwidth $b$ on the torus: $F_b(t) = C_b \exp(-d \times (0,t)^2/b^2)$, where $d(t,t')$ is the flat Riemannian distance on the torus and $C_b$ is a normalizing constant. We then compare the different thresholds obtained by the methods proposed in this work for varying values of $b$. Remember that the algorithms considered here have no prior knowledge on the specific form of the function $F_b$ and would work in other, more complex, dependency contexts.

We consider the two-sided case only. In all of the simulations to come, we fix the following parameters: the dimension is $K = 128^2 = 16{,}384$, the number of data points per sample is $n = 1000$ (hence significantly smaller than $K$) and the width $b$ takes even values in the range $[0, 40]$ ($b = 0$ is white noise; see the left-hand side of Figure 1 for an example of noise realization when $b = 18$). The target test level is $\alpha = 0.05$. We report the (empirical) expectation of each threshold over 250 draws of **Y**.

For computation of the thresholds (9), (11) and (12), we have to choose some parameters $\delta \in (0,1)$ and (for the two latter ones) $\alpha_0 < \alpha$. In each case, these parameters establish a trade-off between a main term and a remainder term; generally speaking, as $n$ grows, one should choose $\delta \to 0$ and $\alpha_0 \to \alpha$ so that the level of the main resampled term tends to the target level $\alpha$. In [1], it was suggested to take $\delta$ of order $1/n$ and $(1 - \frac{\alpha_0}{\alpha})$ of order $n^{-\gamma}$ for some $\gamma > 0$ to ensure that the remainder terms are indeed of lower order, but there is no exact recommendation for fixed $n$. In all of the simulations below, we decided to fix $\delta = (1 - \frac{\alpha_0}{\alpha}) = 0.1$, without particularly trying to optimize this choice. When varying these parameter values, we noticed that the results were not overly sensitive to them. Finally, for all of the thresholds, the resampling quantities (quantiles or expectations) are estimated by Monte Carlo with 1000 draws (but we disregarded the additional terms proposed in [1] to account for the Monte Carlo random error).

4.2. *Simulations with unspecified alternative: Single-step, translation invariant procedures.* We first study the performance of the multiple testing procedures which have a translation invariant threshold (TI), that is, the single-step procedures using thresholds (7), (8), (9), (11) and (12), denoted, respectively, by "bonf," "conc," "conc ∧ bonf," "quant + bonf" and "quant + conc." Their distributions do not depend on the true mean vector $\mu$ chosen to generate data and we have fixed $\mu = 0$ without loss of generality. Provided that the FWER constraint is satisfied, procedures with a smaller threshold are less conservative and hence more powerful.

In Figure 1, we report the (averaged) values of each threshold. In this figure, we did not include standard deviations: they are quite low, of the order of $10^{-3}$, although it is worth noting that the quantile threshold has a standard deviation roughly twice as large as the concentration threshold. For comparison, we also included an estimation of the true quantile, that is, the $1 - \alpha$ quantile of the distribution of $\sup_{k \in \mathcal{H}} |\overline{\mathbf{Y}}_k - \mu_k|$ (more precisely, an empirical quantile over 1000 samples), denoted by "ideal." The exact threshold corresponding to $K = 1$ (test of a single coordinate Gaussian mean) is also included for comparison and is denoted by "single." In the context of this experiment, we also computed the threshold (10), that is, the raw symmetrized quantile obtained after empirical recentering of the data (for which no nonasymptotic theoretical results are available). This threshold was not reported in the plots because it turns out to be so close to the true quantile that they are almost indistinguishable.

The overall conclusion of this first experiment is that the different thresholds proposed in this work are relevant: they improve on the Bonferroni threshold, provided the vector has strong enough correlations. As expected, the quantile approach appears to lead to tighter thresholds. (This might, however, not always be the case for smaller sample sizes because of the additional term.) One remaining advantage of the concentration approach is that the compound threshold (9) falls back on the Bonferroni threshold when needed, at the cost of a minimal threshold increase. Finally, the remainder terms introduced by the theory in the centered quantile thresholds appear overestimated since the raw resampled quantile is, in fact, extremely close to the true quantile.

4.3. *Simulations with a specific alternative.* We consider the experiment of the previous section, with the following choice for the vector of true means:

(16) $\quad \forall (i,j) \in \{0, \ldots, 127\}^2 \qquad \mu_{(i,j)} = \frac{(64-j)_+}{64} \times 20 t_{\alpha,\mathrm{Bonf}}(\mathcal{H}).$

In this situation, half of the null hypotheses are true, while the nonzero means are increasing linearly from $(5/16)t_{\alpha,\mathrm{Bonf}}(\mathcal{H})$ to $20 t_{\alpha,\mathrm{Bonf}}(\mathcal{H})$. The thresholds obtained are displayed in Figure 2, along with the averaged power of the corresponding procedures, defined as the expected proportion of signal correctly detected (i.e., averaged proportion of rejections among the false null hypotheses).

In this experiment, we concentrated on the quantile-based thresholds. We picked the threshold (12) "quant + bonf" as a representative of the centered quantile approach and its step-down counterpart, denoted "s.d. quant + bonf." We compare these to the uncentered quantile threshold (13), denoted "quant.uncent," and its step-down version, "s.d. quant.uncent." Bonferroni's threshold and its step-down version "holm" are included for comparison. The threshold denoted "ideal" is now derived from the $1-\alpha$ quantile of the distribution of $\sup_{k, \mu_k = 0} |\overline{\mathbf{Y}}_k|$ and corresponds to the optimal threshold for FWER control.

The results of the experiment can be summarized as follows:

- The single-step centered quantile procedure "quant + bonf" outperforms Holm's procedure provided the coordinates of the vector are sufficiently correlated. Its step-down version "s.d. quant + bonf" performs even better, although the difference is not huge.
- The single-step procedure based on the uncentered quantile "quant.uncent" has the worst performance, confirming the qualitative analysis following Theorem 2.3.

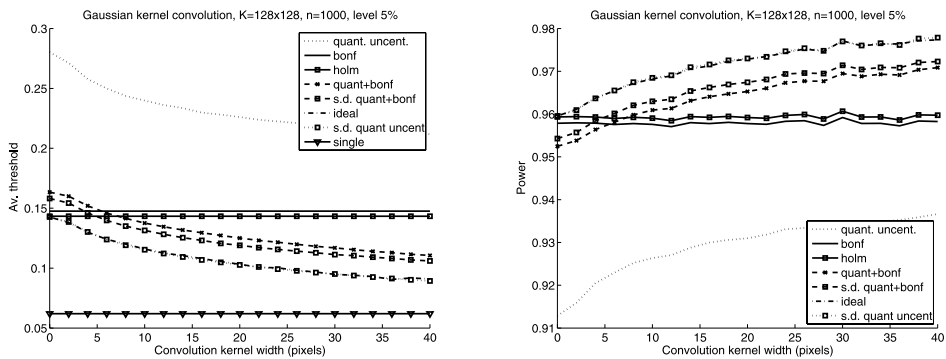

Fig. 2. *Multiple testing problem with $\mu$ defined by (16) for different approaches; see text. Left: average thresholds. Right: average power.*

- The step-down procedure based on the uncentered quantile "s.d. quant.uncent" seems, on the other hand, to be the most accurate among the procedures considered here, also in accordance with the qualitative analysis following Theorem 3.2.

The latter point must be balanced with computation time considerations. When $K$ and $n$ are large, the step-down algorithm for the uncentered quantile takes longer to compute because of its iterative nature, while the single-step centered quantile procedure "quant + bonf" provides a relatively good accuracy without iterating. This brings us to the next point.

4.4. *Hybrid approach.* We show here, with a specific simulation study, that the hybrid approach proposed in Algorithm 3.3 results in a speed/accuracy trade-off which is particularly noticeable when the mean values take on a large range.

Consider the same simulation framework as above, except that the bandwidth $b$ is now fixed at 30, the size of the sample is $n = 100$ and the means are given as follows: $\forall (i,j) \in \{0, \ldots, 127\}^2$, $\mu_{(i,j)} = f(i + 128j)$, where

$$
\begin{aligned}
\forall k \in \{0, \ldots, 128^2/2\} \quad f(k) &= 0.5 t_{\alpha, \mathrm{Bonf}}(\mathcal{H}) \\
&\quad \times \exp\left( \frac{128^2/2 - k}{128^2/2} \frac{r}{10} \log(10) \right)
\end{aligned}
\tag{17}
$$

and $f(k) = 0$ for the other values of $k$. In this situation, the $128^2/2$ nonzero means are decreasing exponentially between $0.5 t_{\alpha, \mathrm{Bonf}}(\mathcal{H}) 10^{r/10}$ and $0.5 t_{\alpha, \mathrm{Bonf}}(\mathcal{H})$, where $r$ is the dynamic range (in dB) of the signal.

In Figure 3, we have computed, for several values of $r$, the average number of iterations for the above step-down procedures, as well as their power when these procedures are stopped early after at most $t$ iterations (such an early stopping is relevant in the case of a strict computation time constraint). We can sum up these results as follows:

- The hybrid approach needs, on average, significantly less iterations to converge when $r \geq 10$.
- Stopping the hybrid approach procedure after only two iterations results in an average power that is virtually indistinguishable from the power obtained without early stopping, uniformly over values of $r$. In contrast, as $r$ increases, more iterations are needed for the step-down quantile uncentered threshold in order to reach full power.

While these results are certainly specific to the particular simulation setup we used, they illustrate that the informal and qualitative analysis presented in Section 3 is correct when the signal (nonzero means) has a wide dynamic range. In particular, the fact that the hybrid approach already gives very satisfactory results after the first two iterations reinforces the interpretation that the first step (using the centered quantile threshold with remainder term) immediately rules out all coordinates with a large SNR, while the second step (using the exact, uncentered quantile) improves the precision once these high-SNR coordinates have been eliminated.

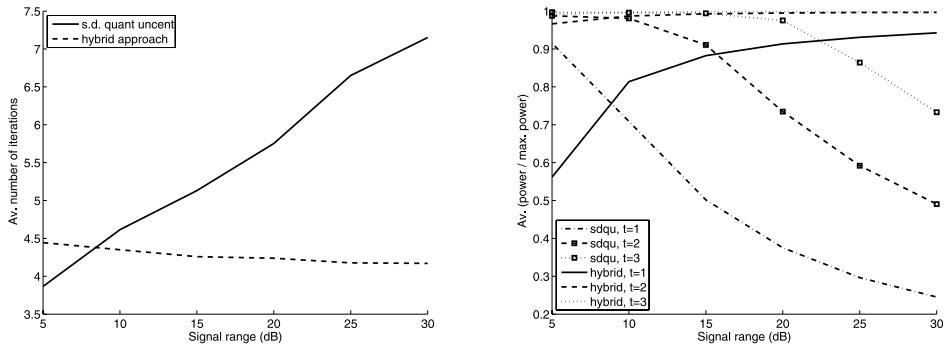

FIG. 3. *Multiple testing problem with $\mu$ corresponding to (17) for the step-down procedure based on the uncentered quantile (sdqu) and the hybrid step-down approach. Left: average number of step-down iterations. Right: average of the ratio of power to maximum power when the step-down is stopped after at most $t$ iterations. Here, the maximum power is taken to be the power of "sdqu" without early stopping. (For the hybrid approach, the first step counts as one iteration.)*

## 5. Discussion and concluding remarks.

5.1. *Discussion: FWER versus FDR in multiple testing.* It can legitimately be asked whether the FWER is an appropriate measure of type I error. The false discovery rate (FDR), introduced in [2] and defined as the average proportion of wrongly rejected hypotheses among all of the rejected hypotheses, appears to have recently become a de facto standard, in particular, in the setting of a large number of hypotheses to test, as we consider here. One reason for the popularity of the FDR is that it is a less strict measure of error than the FWER and, to this extent, FDR-controlled procedures reject more hypotheses than FWER-controlled ones. We give two reasons why the FWER is still a quantity of interest to investigate. First, the FDR is not always relevant, in particular, for neuroimaging data. Indeed, in this context, the signal is often strong over some well-known large areas of the brain (e.g., the motor and visual cortices). Therefore, if, for instance, 95 percent of the detected locations belong to these well-known areas, FDR control (at level 5%) does not provide evidence for any new true discovery. On the contrary, FWER control is more conservative, but each detected location outside these well-known areas is a new true discovery with high probability. Second, assuming that the FDR or a related quantity is nevertheless the end goal, it can be very useful to consider a two-step procedure, where the first step consists of a FWER-controlled multiple test. Namely, this first step can be used as a mean to estimate the FDR or the FDP (false discovery proportion) of another procedure used in the second step and thus to fine-tune the parameters of this second step for the desired goal. This approach has been advocated, for example, in [3, 4] for finding FDR controlling procedures adaptive to the proportion of true nulls and in [12] to find specific regions in random fields, also with application to neuroimaging data.

5.2. *Conclusion.* In this work, the main point was to introduce multiple testing procedures based on resampling thresholds (9), (11) and (12) coming from nonasymptotic confidence regions constructed in [1]. These confidence regions have theoretical control of the confidence level for any $n$, so the FWER of the corresponding multiple testing procedures is also controlled nonasymptotically. This issue is important in practice because the sample size is often much smaller than the number of tests to perform ($K \gg n$). Nevertheless, as the simulations of Section 4 suggest, remainder terms in the thresholds—precisely introduced to deal with this nonasymptotic setting—are overestimated by the theory and could probably be improved.

Even in the presence of these corrective terms, we showed through experiments that these thresholds are able to capture the unknown dependency structure of the data and significantly outperform Holm's procedure when this dependency is strong enough. In comparison to exact randomization tests (based on an uncentered quantile), which also provide nonasymptotic level control, we argued that the empirical centering operation before random sign reversal results in translation invariant thresholds. These thresholds are, for this reason, unaffected by the unknown signal and thus already relevant for testing in the first iteration of the step-down algorithm. The method also applies to one-sided testing problems, where the uncentered approach is not theoretically justified as far as we know. Finally, the hybrid algorithm can approach the accuracy of the uncentered step-down threshold (which does not require corrective terms) while initially taking advantage of the centered threshold, resulting in a faster computation.

For practical purposes, it is certainly tempting to recommend using a (step-down) procedure based on the raw, unmodified centered quantile without remainder terms (10): this would correspond to the principle of traditional resampling. To this extent, and to rephrase the discussion in [1], nonasymptotic theoretical results can also be understood from an asymptotic point of view, justifying the use of resampling (in a specific setting— Gaussian variables, test for the mean, Rademacher weights) for a regime that is not usually covered by traditional asymptotics (i.e., dimension $K_n$ increasing with $n$).

**6. Proof of Proposition 3.4.** First, note that $q_{\alpha_0}(\mathbf{Y}, \mathcal{H}_0) \leq q_{\alpha_0}(\mathbf{Y} - \mu, \mathcal{H})$. From the proof of Theorem 3.2 in [1], with probability greater than $1 - (\alpha - \alpha_0)$, we have

$$q_{\alpha_0}(\mathbf{Y} - \mu, \mathcal{H}) \leq t_{\alpha, \text{quant}+\text{Bonf}}(\mathbf{Y}, \mathcal{H}).$$

Take $\mathbf{Y}$ in the event where the above inequality holds. If the global procedure rejects at least one true null hypothesis, let $j_0$ denote the first time that this occurs ($j_0 = 0$ if it is in the first step). There are two cases:

- if $j_0 = 0$, then $\sup_{k \in \mathcal{H}_0} |\overline{\mathbf{Y}}_k| \geq t_{\alpha, \text{quant}+\text{Bonf}}(\mathbf{Y}, \mathcal{H}) \geq q_{\alpha_0}(\mathbf{Y} - \mu, \mathcal{H}) \geq q_{\alpha_0}(\mathbf{Y}, \mathcal{H}_0)$;
- if $j_0 \geq 1$, then $\sup_{k \in \mathcal{H}_0} |\overline{\mathbf{Y}}_k| \geq t_{\alpha_0, \text{quant.uncent}}(\mathbf{Y}, \mathcal{C}_{j_0-1})$ and $\mathcal{H}_0 \subset \mathcal{C}_{j_0-1}$ (from the definition of $j_0$) so that $\sup_{k \in \mathcal{H}_0} |\overline{\mathbf{Y}}_k| \geq t_{\alpha_0, \text{quant.uncent}}(\mathbf{Y}, \mathcal{H}_0) = q_{\alpha_0}(\mathbf{Y}, \mathcal{H}_0)$.

In both cases, $\sup_{k \in \mathcal{H}_0} |\overline{\mathbf{Y}}_k| \geq q_{\alpha_0}(\mathbf{Y}, \mathcal{H}_0)$, which occurs with probability smaller than $\alpha_0$.

**Acknowledgments.** The first author's research was mostly carried out at University Paris-Sud (Laboratoire de Mathematiques, CNRS UMR 8628). The second author's research was partially carried out while holding an invited position at the Statistics Department of the University of Chicago, which is warmly acknowledged. The third author's research was mostly carried out at the French institute INRA-Jouy and at the Free University of Amsterdam. We would like to thank the two referees and the Associate Editor for their insights which led, in particular, to a more rational organization of the paper.

## REFERENCES

[1] ARLOT, S., BLANCHARD, G. and ROQUAIN, É. (2010). Some nonasymptotic results on resampling in high dimension, I: Confidence regions. *Ann. Statist.* **38** 51–82.
[2] BENJAMINI, Y. and HOCHBERG, Y. (1995). Controlling the false discovery rate: A practical and powerful approach to multiple testing. *J. Roy. Statist. Soc. Ser. B* **57** 289–300. MR1325392
[3] BLANCHARD, G. and ROQUAIN, É. (2008). Two simple sufficient conditions for FDR control. *Electron. J. Stat.* **2** 963–992. MR2448601
[4] BLANCHARD, G. and ROQUAIN, É. (2009). Adaptive FDR control under independence and dependence. *J. Mach. Learn. Res.* To appear.
[5] DARVAS, F., RAUTIAINEN, M., PANTAZIS, D., BAILLET, S., BENALI, H., MOSHER, J., GARNERO, L. and LEAHY, R. (2005). Investigations of dipole localization accuracy in MEG using the bootstrap. *NeuroImage* **25** 355–368.
[6] DUROT, C. and ROZENHOLC, Y. (2006). An adaptive test for zero mean. *Math. Methods Statist.* **15** 26–60. MR2225429
[7] FISHER, R. A. (1935). *The Design of Experiments.* Oliver and Boyd, Edinburgh.
[8] GE, Y., DUDOIT, S. and SPEED, T. P. (2003). Resampling-based multiple testing for microarray data analysis. *Test* **12** 1–77. MR1993286

[9] HOLM, S. (1979). A simple sequentially rejective multiple test procedure. *Scand. J. Statist.* **6** 65–70. MR0538597
[10] JERBI, K., LACHAUX, J.-P., N'DIAYE, K., PANTAZIS, D., LEAHY, R. M., GARNERO, L. and BAILLET, S. (2007). Coherent neural representation of hand speed in humans revealed by MEG imaging. *PNAS* **104** 7676–7681.
[11] PANTAZIS, D., NICHOLS, T. E., BAILLET, S. and LEAHY, R. M. (2005). A comparison of random field theory and permutation methods for statistical analysis of MEG data. *NeuroImage* **25** 383–394.
[12] PACIFICO, M. P., GENOVESE, I., VERDINELLI, I. and WASSERMAN, L. (2004). False discovery control for random fields. *J. Amer. Statist. Assoc.* **99** 1002–1014. MR2109490
[13] POLLARD, K. S. and VAN DER LAAN, M. J. (2004). Choice of a null distribution in resampling-based multiple testing. *J. Statist. Plann. Inference* **125** 85–100. MR2086890
[14] ROMANO, J. P. (1989). Bootstrap and randomization tests of some nonparametric hypotheses. *Ann. Statist.* **17** 141–159. MR0981441
[15] ROMANO, J. P. (1990). On the behavior of randomization tests without a group invariance assumption. *J. Amer. Statist. Assoc.* **85** 686–692. MR1138350
[16] ROMANO, J. P. and WOLF, M. (2005). Exact and approximate stepdown methods for multiple hypothesis testing. *J. Amer. Statist. Assoc.* **100** 94–108. MR2156821
[17] ROMANO, J. P. and WOLF, M. (2007). Control of generalized error rates in multiple testing. *Ann. Statist.* **35** 1378–1408. MR2351090
[18] WESTFALL, P. H. and YOUNG, S. S. (1993). *Resampling-Based Multiple Testing: Examples and Methods for P-Value Adjustment.* Wiley, New York.
[19] YEKUTIELI, D. and BENJAMINI, Y. (1999). Resampling-based false discovery rate controlling multiple test procedures for correlated test statistics. *J. Statist. Plann. Inference* **82** 171–196. MR1736442

S. ARLOT
CNRS: WILLOW PROJECT-TEAM
LABORATOIRE D'INFORMATIQUE
  DE L'ECOLE NORMALE SUPERIEURE
(CNRS/ENS/INRIA UMR 8548)
23 AVENUE D'ITALIE, CS 81321
75214 PARIS CEDEX 13
FRANCE
E-MAIL: sylvain.arlot@ens.fr

G. BLANCHARD
WEIERSTRASS INSTITUTE
  FOR APPLIED STOCHASTICS
  AND ANALYSIS
MOHRENSTRASSE 39, 10117 BERLIN
GERMANY
E-MAIL: blanchar@wias-berlin.de

E. ROQUAIN
UPMC UNIVERSITY OF PARIS 6
UMR 7599, LPMA
4, PLACE JUSSIEU
75252 PARIS CEDEX 05
FRANCE
E-MAIL: etienne.roquain@upmc.fr